\documentclass[10pt]{amsart}
\usepackage{amsmath, amssymb, amsthm, a4wide}
\usepackage{geometry}
\usepackage{fancyhdr}
\usepackage{enumitem}
\newtheorem{theorem}{Theorem}[section]
\newtheorem{lemma}[theorem]{Lemma}
\newtheorem{definition}[theorem]{Definition}
\newtheorem{proposition}[theorem]{Proposition}

\newtheorem{remark}[theorem]{Remark}

\numberwithin{equation}{section}
\DeclareMathOperator*{\esssup}{\textnormal{ess\,sup}}
\DeclareMathOperator*{\essinf}{\textnormal{ess\,inf}}
\DeclareMathOperator*{\dist}{\textnormal{dist}}

\begin{document}

\title[Quadratic Second Order BSDEs]{A New Result for Second Order BSDEs \\ with Quadratic Growth and its applications}
\author{Yiqing LIN}
\address{Institut de Recherche Math\'{e}matique de Rennes, \newline\indent Universit\'{e} de Rennes 1, \newline\indent 35042 Rennes Cedex, France}
\email{yiqing.lin@univ-rennes1.fr}

\date{November 30th, 2012}



\begin{abstract}
In this paper, we study a class of second order backward stochastic differential equations (2BSDEs) with quadratic growth in coefficients. We first establish solvability for such 2BSDEs and then give their applications to robust utility maximization problems.
\end{abstract}
\maketitle
\section{Introduction}
\noindent Typically, nonlinear backward stochastic differential equations (BSDEs) are defined on a Wiener probability space $(\Omega, \mathcal{F}, \mathbb{P})$ and of the following type:
\begin{equation}\label{BSDE}
Y_t= \xi+\int^T_t g(s, Y_s, Z_s)ds-\int^T_tZ_sdB_s,\ 0\leq t \leq T, 
\end{equation}
where $B$ is a Brownian motion, $\mathcal{F}$ is the $\mathbb{P}$-augmented natural filtration generated by $B$,
$g$ is a nonlinear generator, $T$ is the terminal time and $\xi\in\mathcal{F}_T$ is the terminal value. A solution to BSDE (\ref{BSDE}) is a couple of processes $(Y, Z)$ adapted to the filtration $\mathcal{F}$. \\[6pt]
\noindent Under a Lipschitz condition on the generator $g$, Pardoux and Peng \cite{PP} first provided the wellposedness of (\ref{BSDE}). Since then, the theory of nonlinear BSDEs has been extensively studied in the past twenty years. Among all the contributions, we only quote the results which is highly related to our present work.\\[6pt] 
A weaker assumption on the generator is that $g$ has a quadratic growth in $z$.
This kind of BSDEs with bounded terminal value condition was first examined by Kobylanski \cite{KO}, who used a weak convergence technique borrowed from PDE literatures to prove the existence and also obtained the uniqueness result under some additional condition on $g$. With the help of contraction mapping principle, Tevzadze \cite{TEV} re-considered this type of BSDEs when the terminal value $\xi$ is small enough in norm. The advantage of the method adopted by Tevzadze \cite{TEV} is its applicability to not only one-dimensional quadratic BSDEs but also to multidimensional ones. Particularly, the restriction on $\xi$ can be loosen when $g$ satisfies some restrictive condition on regularity. Briand and Hu \cite{BH1, BH2} extended the existence result for (\ref{BSDE}) to the case that $\xi$ is not uniformly bounded and provided the uniqueness result when $g$ is convex. Besides, Morlais \cite{M} considered some similar type of BSDEs driven by continuous martingales. \\[6pt]
Motivated by expected utility theory, Peng \cite{P1} defined a so-called $g$-expectation $\mathcal{E}^g[\xi]:=Y_0$ on $\mathcal{F}_T$ via nonlinear BSDEs with Lipschitz generator. Also, a conditional expectation can be consistently defined: $\mathcal{E}^g[\xi|\mathcal{F}_t]:=Y_t$, under which the solution $Y$ of the BSDE with the generator $g$ is a $g$-martingale. As the counterparts in the classical framework under a linear expectation, Peng \cite{P2} gave the notion of $g$-supermartingle ($g$-submartingle) and established the nonlinear Doob-Meyer type decomposition theorem. Subsequently, Chen and Peng \cite{CP} proved the downcrossing inequality for $g$-martingales. For the case that $g$ is allowed to have a quadratic growth in $z$, similar results can be found in Ma and Yao \cite{MY}.\\[6pt]
Recently, Soner et al. \cite{STZ3} established a framework of ``quasi-sure'' stochastic analysis under a non-dominated class of probability measures. This provided a new approach for Soner et al. \cite{STZ1, STZ2} to re-consider the wellposedness of second order BSDEs (2BSDEs) introduced by Cheridito et al. \cite{CST}. The key idea in Soner et al. \cite{STZ2} is to reinforce a condition that the following 2BSDE holds true $\mathcal{P}_H$-quasi-surely, i.e., $\mathbb{P}$-a.s. for all $\mathbb{P} \in\mathcal{P}_H$, which is a class of mutually singular probability measures (cf. Definition \ref{qs}):
\begin{equation}\label{2BSDEintro}
Y_t=\xi+\int^1_t\hat{F}_s(Y_s, Z_s)ds-\int^1_tZ_sdB_s+K_1-K_t,\ 0\leq t\leq 1.
\end{equation}
Under a uniformly Lipschitz condition on the generator $\hat{F}$, Soner et al. \cite{STZ2} provided a complete wellposedness result for the 2BSDE (\ref{2BSDEintro}). In this pioneering work, a representation theorem of the solution $Y$ is established and thus, the uniqueness is a straightforward corollary. For the existence, a process $Y$ is pathwisely constructed and verified as a $\hat{F}$-supermartingale under each $\mathbb{P}\in\mathcal{P}_H$. Applying the nonlinear Doob-Meyer decomposition theorem, the right-hand side of (\ref{2BSDEintro}) 
comes out, where
$K$ is a (family of) non-decreasing process(es) that satisfies the minimum condition (cf. Definition \ref{def2bsde}).
Moreover, both Cheridito et al. \cite{CST} and Soner et al. \cite{STZ2} explained the connection between the Markov 2BSDEs and a large class of fully nonlinear PDEs, which was one of the motivations initiate this 2BSDEs topic.\\[6pt]
Meanwhile, Peng \cite{P3, P4} independently introduced another framework (so-called $G$-framework) of a time consistent nonlinear expectation $\mathbb{E}_G[\cdot]$, in which a new type of Brownian motion was constructed and the related It\^o type stochastic calculus was established. By explicit constructions,   
Denis et al. \cite{DHP} showed that
$G$-expectation is in fact an upper expectation related to a non-dominated family $\mathcal{P}_G$ 
that consists of some probability measures similar to the elements in $\mathcal{P}_H$. In this regards, the $G$-framework is highly related to the 2BSDE one. 
Adopted the idea in Denis and Martini \cite{DM}, Denis et al. \cite{DHP} defined a Choquet capacity $\bar{C}(\cdot)$ on $(\Omega, \mathcal{B}(\Omega))$ as follows:
$$
\bar{C}(A):=\sup_{\mathbb{P}\in\mathcal{P}_G}\mathbb{P}(A),\ A\subset\mathcal{B}(\Omega),
$$
and then they gave the 
the notion of ``quasi-surely'' in a standard capacity-related vocabulary: a property holds true quasi-surely if and only if it holds outside a polar set, i.e., outside a set $A\subset\Omega$ that satisfies $\bar{C}(A)=0$. We notice that this notion is a little bit stronger than the corresponding one in the 2BSDE framework, so that it yields another type of ``quasi-sure'' stochastic analysis. In this $G$-framework, Hu et al. \cite{HJPS} have worked on nonlinear BSDEs driven by $G$-brownian motion (GBSDEs), which is of the same form as (\ref{2BSDEintro}) but holds in the stronger ``quasi-sure'' sense. In that paper, the solution is an aggregated triple $(Y, Z, K)$ which quasi-surely solves (\ref{2BSDEintro}), where $-K$ is a decreasing $G$-martingale that comes from the $G$-martingale decomposition. 
To ensure that (\ref{2BSDEintro}) is well defined in $G$-framework, an additional condition to the Lipschitz one is imposed on the regularity of the generator (cf. (H1) in Hu et al. \cite{HJPS}). This cost is intelligible since the definition of $G$-stochastic integrals is under a stronger norm induced by $\mathbb{E}_G[\cdot]$ and it makes the space of admissible integrands smaller than the classical one. \\[6pt]
Following the works of Soner et al. \cite{STZ1, STZ2, STZ3}, Possamai and Zhou \cite{PZ} generalized the existence and uniqueness results for the 2BSDE whose generator has a quadratic growth. Based on the previous one of Tevzadze \cite{TEV} for quadratic BSDEs, this work requires some additional condition, either on the terminal value or on the regularity of the generator. Our aim of the present paper is to remove these conditions, that is, to redo the job of Possamai and Zhou \cite{PZ} under some weaker assumptions of the type similar to that in Kobylanski \cite{KO} and Morlais \cite{M}.\\[6pt]
In the classical framework, the quadratic BSDE is a powerful technique to deal with the utility maximization problems. El Karoui and Rouge \cite{ER} computed the value function of an exponential utility maximization problem when the strategies are confined to a convex cone, and they found that its dual problem is related to a quadratic BSDE. In contrast to this, Hu et al. \cite{HIM} and Morlais \cite{M} directly treated the primal problem rather than the dual one and obtained an similar result without the convex condition on the constrain set. The value function was characterized by also a solution of a quadratic BSDE.\\[6pt]
Corresponding to thses works above, Matoussi et al. \cite{MPZ} found that a robust utility maximization problem with non-dominated models can be solved via the 2BSDE technique. This kind of problem was first consider by Denis and Kervarec \cite{DK} under a weakly compact class of probability measures. Just because of this weakly compact assumption, one can find a least favorable probability in this class and work under this probability to find an optimal strategy similarly to how we solve the classical problem under a single probability. With the help of 2BSDEs, Matoussi et al. \cite{MPZ} solved this problem globally and characterized the value function by using a solution of a 2BSDE. This method does not require that the class is weakly compact. However, the result in Matoussi et al. \cite{MPZ} has some limitations: for example, when the utility function is exponential, they are able to solve only the case that $\xi$ is small enough or the border of the constraint domain satisfies an extra regularity condition. This limitations is derived from the theory of quadratic 2BSDEs in Possamai and Zhou \cite{PZ}. Since we shall remove these extra conditions adopted by Possamai and Zhou \cite{PZ}, we can have a better result on solving this robust utility maximization problem.\\[6pt]
This paper is organized as follows: Section 2 includes preliminaries for 2BSDEs theory. Section 3 introduce a priori estimates, a representation theorem and the uniqueness result for 2BSDEs with quadratic growth. Section 4 studies the existence of solutions while Section 5 is the applications of quadratic 2BSDEs to robust maximization problems.
\section{Preliminaries}
\noindent The aim of this section is to list some basic definitions for 2BSDEs introduced by Soner et al. \cite{STZ1, STZ2, STZ3} and Possamai and Zhou \cite{PZ}. The reader interested in a more detailed description of these notation is referred to these papers listed above.
\subsection{The class of probability measures}
Let $\Omega:=\{\omega: \omega\in \mathcal{C}([0, 1], \mathbb{R}^d), \omega_0=0\}$ be the canonical space equipped with the uniform norm $||\omega||^\infty_1:=\sup_{0\leq t\leq 1}|\omega_t|$, $B$ the canonical process, $\mathcal{F}$ the filtration generated by $B$, $\mathcal{F}^+$  the right limit of $\mathcal{F}$.\\[6pt]
We call $\mathbb{P}$ a local martingale measure if under which  
the canonical process $B$ is a local martingale.
By Karandikar \cite{KAR}, the quadratic variation process of $B$ and its density can be defined universally, such that under each local martingale measure $\mathbb{P}$:
%
$$\langle B \rangle:=B^2_t-2\int^t_0B_sdB_s\ {\rm and}\
\hat{a}_t:= \overline{\lim_{\varepsilon\downarrow{0}}}\frac{1}{\varepsilon}(\langle B \rangle_t-\langle B \rangle_{t-\varepsilon}),\ 0\leq t\leq 1,\ \mathbb{P}-a.s..$$
Adapting to Soner et al. \cite{STZ3}, we denote $\overline{\mathcal{P}}_W$ the collection of all local martingale measures $\mathbb{P}$ such that $\langle B\rangle_t$ is absolutely continuous in $t$ and $\hat{a}$ takes values in $\mathbb{S}_d^{>0}$, $\mathbb{P}$-a.s.. It is easy to verify that the following stochastic integral defines a $\mathbb{P}$-Brownian motion:
$$W^\mathbb{P}_t:=\int^t_0\hat{a}_s^{{-1/2}}dB_s,\ 0\leq t\leq 1.$$
We define a subclass of $\overline{\mathcal{P}}_W$ that consists of the probability measures induced by the strong formulation (cf. Lemma 8.1 in Soner et al. \cite{STZ3}):
$$\overline{\mathcal{P}}_S
:=\{\mathbb{P}\in\overline{\mathcal{P}}_W
:\overline{\mathcal{F}^{W^\mathbb{P}}}^\mathbb{P}=\overline{\mathcal{F}}^\mathbb{P}
\},$$
where $\overline{\mathcal{F}}^\mathbb{P}$ ($\overline{\mathcal{F}^{W^\mathbb{P}}}^\mathbb{P}$, respectively) is the $\mathbb{P}$-augmentation of the filtration generated by $B$ ($W^\mathbb{P}$, respectively).
\subsection{The nonlinear generator}
We consider a mapping $H_t(\omega, y, z, \eta): [0,1]\times\Omega\times\mathbb{R}\times\mathbb{R}^d\times D_H\rightarrow \mathbb{R}$ and its Fenchel-Legendre conjugate with respect to $\eta$:
$$
F_t(\omega, y, z, a):=\sup_{\eta\in D_H}\bigg\{\frac{1}{2}{\rm tr}(a\eta)-H_t(\omega, y, z, \eta)\bigg\},\ a\in \mathbb{S}_d^{>0}.
$$
where $D_H\subset\mathbb{R}^{d\times d}$ a given subset that contains $0$. For simplicity of notation, we note
$$\hat{F}_t(y, z):= F_t(y, z, \hat{a}_t)\ {\rm and}\ F^0_t:=\hat{F}_t(0,0),$$
and we denote by $D_{F_t(y, z)}$ the domain of $F$ in $a$ for a fixed $(t,\omega, y, z)$. In accordance with the settings previous literatures, we assume the following assumptions on $F$, which is needed for the ``quasi-sure'' technique:\\[6pt]
\textbf{(A1)} \textit{$D_{F_t(y,z)}=D_{F_t}$ is independent of $(\omega, y, z)$;}\\[6pt]
\textbf{(A2)} \textit{$F$ is $\mathcal{F}$-progressively measurable and uniformly continuous in $\omega$.}
\subsection{The spaces and the norms}
For the wellposedness of 2BSDEs, we consider a restrictive subclass $\mathcal {P}_H\subset\mathcal{P}_S$ defined as follows:
\begin{definition}\label{qs}
Let $\mathcal{P}_H$ denote the collection of all those $\mathbb{P}\in \overline{\mathcal{P}}_S$ such that $$\underline{a}^\mathbb{P}\leq\hat{a}_t\leq\overline{a}^\mathbb{P}\ {(usual\ partial\ ordering\ on}\ \mathbb{S}^{>0}_d {)\ and}\ \hat{a}_t\in D_{F_t},\ \lambda\times \mathbb{P}-a.e.,
$$
for some $\underline{a}^\mathbb{P}$, $\overline{a}^\mathbb{P}\in\mathbb{S}_d^{>0}$ and all $(y,z)\in\mathbb{R}\times\mathbb{R}^d$.
\end{definition}
\begin{remark} Soner et al. \cite{STZ2} mentioned that the bounds $\underline{a}^\mathbb{P}$ and  $\overline{a}^\mathbb{P}$ may vary in $\mathbb{P}$. Thanks to the quadratic growth assumption on $F$, i.e., (A3) in the sequel, $\hat{F}^0_t$ is bounded so that $\mathcal{P}_H$ is not empty in our case (cf. Remark 2.5 in Possamai and Zhou \cite{PZ}).
\end{remark}
\begin{definition}
We say that a property holds $\mathcal{P}_H$-quasi-surely ($\mathcal{P}_H$-q.s.) if it holds $\mathbb{P}$-a.s. for all $\mathbb{P}\in\mathcal{P}_H$.
\end{definition}
\noindent For each $p\geq1$, $L^p_H$ denotes the space of all $\mathcal{F}_1$-measurable scalar random variable $\xi$ that satisfies
$$||\xi||_{L^p_H}:=\sup_{\mathbb{P}\in \mathcal{P}_H}\mathbb{E}^{\mathbb{P}}[|\xi|^p]<+\infty.$$
Letting $p\rightarrow+\infty$, we denote by $L^\infty_H$ the space of all $\mathbb{P}_H$-q.s. bounded random variable $\xi$ with 
$$||\xi||_{L^\infty_H}:=\sup_{\mathbb{P}\in\mathcal{P}_H}||\xi||_{L^\infty(\mathbb{P})}<+\infty.
$$
\noindent Let $\mathbb{D}^\infty_H$ denote the space of all $\mathbb{R}$-valued $\mathcal{F}^+$-progressively measurable process $Y$ that satisfies
$$\mathcal{P}_H-q.s.\ c\grave{a}dl\grave{a}g\ {\rm and}\ ||Y||_{D^\infty_H}:=\sup_{0\leq t\leq1}||Y_t||_{L^\infty_H}<+\infty,$$
and $\mathbb{H}^2_H$ denotes the space of all $\mathbb{R}^d$-valued $\mathcal{F}^+$-progressively measurable process $Z$ that satisfies
$$||Z||^2_{\mathbb{H}^2_H}:=\sup_{\mathbb{P}\in\mathcal{P}_H}\mathbb{E}^\mathbb{P}\bigg[\int^1_0|\hat{a}^{1/2}_tZ_t|^2dt\bigg]<+\infty.$$
\begin{remark}
We emphasize that the monotone convergence theorem no long holds true on each space listed above in this framework, i.e. that the monotone $\mathcal{P}_H$-q.s. convergence yields the convergence in norm may fail.
As stated in section 4 of Possamai and Zhou \cite{PZ}, this is one of the main difficulties to prove the existence of quadratic 2BSDEs by global approximation.
\end{remark}
\noindent With a little abuse of notation, we introduce the notion of $BMO(\mathcal{P}_H)$-martingale and its generator, which is an extension of the classical one.  For the convenience of notation, $H$ can refer to either a single process or a family of non-aggregated processes $\{H^\mathbb{P}\}_{\mathbb{P}\in\mathcal{P}_H}$ in the definition and lemmas below.
\begin{definition}\label{bmodef}
We call $H$ a $BMO(\mathcal{P}_H)$-martingale if for each $\mathbb{P}\in\mathcal{P}_H$, $H^\mathbb{P}$ is a $\mathbb{P}$-square integrable martingale and 
$$
||H||^2_{BMO_2(\mathcal{P}_H)}:=\sup_{\mathbb{P}\in\mathcal{P}_H}\sup_{\tau\in\mathcal{T}^1_0}
||\mathbb{E}^\mathbb{P}_\tau[\langle H^\mathbb{P}\rangle_1-\langle H^\mathbb{P}\rangle_\tau]||_{L^\infty(\mathbb{P})}<+\infty,
$$
where $\mathcal{T}^1_0$ is the collection of all $\mathcal{F}$-stopping times $\tau$ that take values in $[0, 1]$.
\end{definition}
\noindent From the definition above, for a fixed $BMO(\mathcal{P}_H)$-martingale $H$, there exists a uniform bound constant $M_H>0$, such that for all $\mathbb{P}\in\mathcal{P}_H$ and $\sigma\in\mathcal{T}^1_0$,
$$||H_{\cdot\wedge \sigma}||^2_{BMO_2(\mathbb{P})}\leq||H||^2_{BMO_2(\mathbb{P})}
\leq M_H.$$
Applying Theorem 2.4 and Theorem 3.1 in Kazamaki \cite{KAZ} under each $\mathbb{P}\in\mathcal{P}_H$, we have the following lemmas:
\begin{lemma}\label{BMOP}
Suppose $H$ is a $BMO(\mathcal{P}_H)$-martingale,
then there exist two constants $r>1$ and $C>0$, such that
$$\sup_{\mathbb{P}\in\mathcal{P}_H}\sup_{0\leq t\leq 1}\mathbb{E}^\mathbb{P}[|
\mathcal{E}(H^\mathbb{P})
_t|^r]\leq C,$$
and for some $q>1$, the following reverse H\"older's inequality holds under each $\mathbb{P}\in\mathcal{P}_H$ with a uniform constant $C_{RH}$: for each $0 \leq t_1 \leq t_2\leq 1$,
$$
\mathbb{E}^\mathbb{P}_{t_1}[\mathcal{E}(H^\mathbb{P})
^{q}_{t_2}]\leq C_{RH} \mathcal{E}(H^\mathbb{P})
^{q}_{t_1},\ \mathbb{P}-a.s.,
$$
where $r$, $C$, $q$ and $C_{RH}$ simply depend on $M_H$.
\end{lemma}
\begin{lemma}\label{bmolemma2}
Suppose $H$ is a $BMO(\mathcal{P}_H)$-martingale,
then there exist a $p>1$ and a $C_E>0$ that simply depend on $M_H$, such that for each $t\in[0, 1]$,
$$\sup_{\mathbb{P}\in\mathcal{P}_H}\sup_{\tau\in{\mathcal{T}^t_0}}\bigg|\bigg|\mathbb{E}^\mathbb{P}_\tau\bigg[\bigg(
\frac{\mathcal{E}(H^\mathbb{P})_\tau}{\mathcal{E}(H^\mathbb{P})_t}\bigg)^\frac{1}{p-1}\bigg]\bigg|\bigg|_{L^\infty(\mathbb{P})}\leq C_E.$$
\end{lemma}
%

\begin{definition}
We call $Z\in\mathbb{H}^2_H$ a $BMO(\mathcal{P}_H)$-martingale generator if 
\begin{align*}
||Z||^2_{\mathbb{H}^2_{BMO(\mathcal{P}_H)}}:&=
\sup_{\mathbb{P}\in\mathcal{P}_H}\bigg|\bigg|\int^\cdot_0Z_tdB_t\bigg|\bigg|^2_{BMO_2(\mathbb{P})}\\
&=\sup_{\mathbb{P}\in\mathcal{P}_H}\sup_{\tau\in\mathcal{T}^1_0}\bigg|\bigg|\mathbb{E}^\mathbb{P}_\tau\bigg[\int^1_\tau|\hat{a}^{1/2}_tZ_t|^2dt\bigg]\bigg|\bigg|_{L^\infty(\mathbb{P})}<+\infty.
\end{align*}
It is evident that if $Z$ is a $BMO(\mathcal{P}_H)$-martingale generator, 
defining for each $\mathbb{P}\in\mathcal{P}_H$,
$$H^\mathbb{P}_t:=\int^t_0 Z_sdB_s,\ 0\leq t\leq 1,$$
then $H$ is a $BMO(\mathcal{P}_H)$-martingale. We denote by $\mathbb{H}^2_{BMO(\mathcal{P}_H)}$ the space of all $BMO(\mathcal{P}_H)$-martingale generators.
\end{definition}
\noindent Applying energy inequality under each $\mathbb{P}\in\mathcal{P}_H$, we have the following lemma:
\begin{lemma}\label{ei}
Suppose $Z\in\mathbb{H}^2_{BMO(\mathcal{P}_H)}$, for each $p\geq 1$, $\mathbb{P}\in\mathcal{P}_H$ and all $\tau\in\mathcal{T}^1_0$,
$$
\mathbb{E}^\mathbb{P}_\tau\bigg[\bigg(\int^1_\tau|\hat{a}^{1/2}_tZ_t|^2dt\bigg)^p\bigg]\leq C_p||Z||^{2p}_{\mathbb{H}^2_{BMO}},\ \mathbb{P}-a.s..
$$
\end{lemma}
\noindent Finally, we denote by $UC_b(\Omega)$ the collection of all bounded and uniformly continuous maps $\xi: \Omega\rightarrow\mathbb{R}$ and denote by $\mathcal{L}^\infty_H$ the closure of $UC_b(\Omega)$ under the norm$||\cdot||_{L^\infty_H}$.
\subsection{Formulation to quadratic 2BSDEs}
We shall consider the 2BSDE of the following form, which is first introduce in Soner et al. \cite{STZ2}:
\begin{equation}\label{key}
Y_t=\xi+\int^1_t\hat{F}_s(Y_s, Z_s)ds-\int^1_tZ_sdB_s+K_1-K_t,\ 0\leq t\leq 1,\ \mathcal{P}_H-q.s..
\end{equation}
In addition to (A1)-(A2), we assume the following conditions on the generator $F$:\\[6pt]
\textbf{(A3)} \textit{$F$ is continuous in $(y, z)$ and has a quadratic growth, i.e. there exists a triple $(\alpha, \beta, \gamma)\in \mathbb{R}^+\times\mathbb{R}^+\times\mathbb{R}^+$, such that for all $(\omega, t, y, z, a)\in\Omega\times[0, 1]\times\mathbb{R}\times\mathbb{R}^d\times D_{F_t}$,}
\begin{equation}\label{gere}
|F_t(\omega, y, z, a)|\leq \alpha+\beta|y|+\frac{\gamma}{2}|a^{1/2}z|^2;
\end{equation}
\textbf{(A4)} \textit{$F$ is uniform Lipschitz in $y$, i.e. there exists a $\mu>0$, such that for all $(\omega, t, y, y', z, \\a)\in\Omega\times[0, 1]\times\mathbb{R}\times\mathbb{R}\times\mathbb{R}^d\times D_{F_t}$,}
$$
|F_t(\omega, y, z, a)-F_t(\omega, y', z, a)|\leq \mu |y-y'|
;$$
\noindent\textbf{(A5)} \textit{$F$ is local Lipschitz in $z$, i.e. 
for each $(\omega, t, y, z, z' a)\in\Omega\times[0, 1]\times\mathbb{R}\times\mathbb{R}^d\times\mathbb{R}^d\times D_{F_t}$,
$$|F_t(\omega, y, z, a)-F_t(\omega, y, z', a)|\leq C (1+|a^{1/2}z|+|a^{1/2}z'|)|a^{1/2}(z-z')|
.$$
}
%
\begin{remark}
We have some comments on these conditions above:
(A3) is a quadratic growth condition for the proof of existence and similar ones for quadratic BSDEs can be found in Kobylanski \cite{KO} and Marlais \cite{M}; (A4) and (A5) are necessary for the proof of uniqueness and analogous conditions were adopted by Hu et al. \cite{HIM} and Morlais \cite{M} for quadratic BSDEs. All these conditions above could be slightly weakened and further discussion will be made in Remark \ref{rem1}.
\end{remark}
\begin{definition}\label{def2bsde}
We say that $(Y, Z)\in\mathbb{D}^\infty_H\times\mathbb{H}^2_H$ is a solution of 2BSDE (\ref{key}) 
if:\\[6pt]
-\ $Y_T=\xi$, $\mathcal{P}_H$-q.s.;\\[6pt]
-\ The process $K^\mathbb{P}$ defined as below: for each $\mathbb{P}\in\mathcal{P}_H$,
\begin{equation}\label{k}
K^\mathbb{P}_t:=Y_0-Y_t-\int^t_0\hat{F}_s(Y_s, Z_s)ds+\int^t_0Z_sdB_s,\ 0\leq t\leq 1,\ \mathbb{P}-a.s.,
\end{equation}
has non-decreasing paths $\mathbb{P}$-a.s.;\\[6pt]
-\ The family $\{K^\mathbb{P}\}_{\mathbb{P}\in\mathcal{P}_H}$ satisfies the minimum condition: for each $\mathbb{P}\in\mathcal{P}_H$,
\begin{equation}\label{min}
K^\mathbb{P}_t=\sideset{}{^\mathbb{P}}\essinf_{\mathbb{P}'\in\mathcal{P}_H(t^+,\mathbb{P})}\mathbb{E}^{\mathbb{P}'}_t[K^{\mathbb{P}'}_T],\ 0\leq t\leq 1,\ \mathbb{P}-a.s..
\end{equation}
Moreover, if the family $\{K^\mathbb{P}\}_{\mathbb{P}\in\mathcal{P}_H}$ can be aggregated into a universal process $K$, we call $(Y, Z, K)$ a solution of 2BSDE (\ref{key}).
\end{definition}
\noindent In the sequel, positive constants $C$ and $M$ vary from line to line.
\section{Representation and uniqueness of solutions to 2BSDEs}
\noindent In this section, we give a representation theorem of solutions to the 2BSDE (\ref{key}) under (A1)-(A5), which is similar to those in Soner et al. \cite{STZ2} and Possamai and Zhou \cite{PZ}. The representation theorem shows the relationship between the solution to the 2BSDE (\ref{key}) and those to quadratic BSDEs with the generator $\hat{F}$ under each $\mathbb{P}\in\mathcal{P}_H$. Also, some a priori estimates to solutions is given which are useful to the proof of the existence.
\subsection{Representation theorem}
Before proceeding the argument, we first introduce a lemma (cf. Lemma 3.1 in Possamai and Zhou \cite{PZ}), the parallel version of which for quadratic BSDEs plays a very important role to show the connection between the boundness of $Y$ and the $BMO$ property of the martingale part $\int^T_\cdot Z_tdB_t$. 
\begin{lemma}\label{bmo}
We assume (A1)-(A3) and $\xi\in L^\infty_H$. If $(Y, Z)\in\mathbb{D}^\infty_H\times\mathbb{H}^2_H$ is a solution to the 2BSDE (\ref{key}), then $Z\in\mathbb{H}^2_{BMO(\mathcal{P}_H)}$ and 
\begin{equation}\label{er1}
||Z||^2_{\mathbb{H}^2_{BMO(\mathcal{P}_H)}}\leq \frac{1}{\gamma^2}e^{4\gamma||Y||_{\mathbb{D}^\infty_H}}(1+2\gamma (\alpha+\beta||Y||_{\mathbb{D}^\infty_H})).
\end{equation}
\end{lemma}
\noindent Consider the following quadratic BSDE under each $\mathbb{P}\in\mathcal{P}_H$:
\begin{equation}\label{st}
y^\mathbb{P}_s=\eta+\int^t_s\hat{F}_u(y^\mathbb{P}_u, z^\mathbb{P}_u)du-\int^t_sz^\mathbb{P}_udB_u,\ 0\leq s\leq t,\ \mathbb{P}-a.s.,
\end{equation}
where $t\in[0, 1]$ and $\eta$ is a $\mathcal{F}_t$-measurable random variable in $L^\infty(\mathbb{P})$. Under (A1)-(A5), the BSDE (\ref{st}) admits a unique solution $(y^\mathbb{P}(t, \eta), z^\mathbb{P}(t, \eta))$ according to Kobylanski \cite{KO} and Morlais \cite{M}.\\[6pt]
Then, we have the following representation theorem for the solution of the 2BSDE (\ref{key}):
\begin{theorem}\label{rep}
Let (A1)-(A5) hold. Assume that $\xi\in L^\infty_H$ and $(Y, Z)\in\mathbb{D}^\infty_H\times\mathbb{H}^2_H$ is a solution of the 2BSDE (\ref{key}). Then, for each $\mathbb{P}\in\mathcal{P}_H$ and all $0\leq t_1\leq t_2\leq 1$,
\begin{equation}\label{er12}
Y_{t_1}=\sideset{}{^\mathbb{P}}\esssup_{\mathbb{P}'\in\mathcal{P}_H(t^+_1, \mathbb{P})}y^{\mathbb{P}'}_{t_1}(t_2, Y_{t_2}),\ \mathbb{P}-a.s.,
\end{equation}
where $$\mathcal{P}_H(t^+_1, \mathbb{P}):=\{
\mathbb{P}'\in {\mathcal{P}_H}: \mathbb{P}'|_{\mathcal{F}^+_{t_1}}=\mathbb{P}|_{\mathcal{F}^+_{t_1}}\}.$$
\end{theorem}
\begin{remark}
Applying Theorem 2.7 (comparison principle) in Morlais \cite{M}, the theorem above also implies a comparison principle for quadratic 2BSDEs.
\end{remark}
\noindent\textbf{Proof:} 
First of all, Lemma \ref{bmo} shows that $Z$ is a $BMO(\mathcal{P}_H)$-martingale generator, then we deduce by the BDG type inequalities, Lemma \ref{ei} and (\ref{er1}) that
 for each $p\geq 1$, $\mathbb{P}\in\mathcal{P}_H$ and all $0\leq t_1\leq t_2\leq 1$, 
\begin{equation}\label{18}\mathbb{E}^\mathbb{P}_{t_1}[(K^\mathbb{P}_{t_2}-K^\mathbb{P}_{t_1})^p]\leq C_p:=C e^{4p\gamma||Y||_{\mathbb{D}^\infty_H}}(1+||Y||^p_{\mathbb{D}^\infty_H}),\ \mathbb{P}-a.s., \end{equation}
Since $\mathbb{P}$ is arbitrary in (\ref{18}), we have
$$\sideset{}{^\mathbb{P}}\esssup_{\mathbb{P}'\in\mathcal{P}_H(t^+_1,\mathbb{P})}\mathbb{E}^{\mathbb{P}'}_{t_1}[(K^{\mathbb{P}'}_{t_2}-K^{\mathbb{P}'}_{t_1})^p]<C_p,\ \mathbb{P}-a.s..$$
\noindent We are now ready to prove that for a fixed $\mathbb{P}\in\mathcal{P}_H$ and all $0\leq t_1\leq t_2\leq 1$,
\begin{equation}\label{zbzd}
Y_{t_1}\leq\sideset{}{^\mathbb{P}}\esssup_{{\mathbb{P}'}\in\mathcal{P}_H(t^+_1, \mathbb{P})}y^{\mathbb{P}'}_{t_1}(t_2, Y_{t_2}),\ \mathbb{P}-a.s..
\end{equation}
Fixing $t_2\in[0, 1]$, for each $\mathbb{P}'\in\mathcal{P}_H(t^+_1, \mathbb{P})$, 
we note
$$
\delta Y^{\mathbb{P}'}:= Y-y^{\mathbb{P}'}(t_2,Y_{t_2})\ {\rm and}\ \delta Z^{\mathbb{P}'}:=Z-z^{\mathbb{P}'}(t_2,Y_{t_2}),
$$
then, for each $t\in[0, t_2]$,
\begin{align*}
\delta Y^{\mathbb{P}'}_t=\int^{t_2}_t\lambda_s\delta Y_sds-\int^{t_2}_t\delta Z_s\hat{a}^{1/2}_s(-\kappa^{\mathbb{P}'}_sds+dW^{\mathbb{P}'}_s)+K^{\mathbb{P}'}_{t_2}-K^{\mathbb{P}'}_t,\ \mathbb{P}'-a.s.,
\end{align*}
where $\lambda$ is a scalar valued process and $\kappa$ is an $\mathbb{R}^d$-valued process defined by $$
\kappa^{\mathbb{P}'}_t=\left\{
\begin{array}{l@{\quad, \quad}l}
\frac{(\hat{F}_t(y^{\mathbb{P}'}_t(t_2, Y_{t_2}), z^{\mathbb{P}'}_t(t_2, Y_{t_2}))-\hat{F}_t(y^{\mathbb{P}'}_t(t_2, Y_{t_2}), Z_t))\hat{a}^{1/2}_t\delta Z_t}{|\hat{a}^{1/2}_t\delta Z_t|^2} &  |\hat{a}^{1/2}_t\delta Z_t| \neq 0;\\
0 & otherwise.
\end{array}\right.
$$
By (A4) and (A5),
we have $||\lambda||_{D^\infty(\mathbb{P}')}\leq \mu$ and $\kappa$ satisfies 
$$
|\kappa^{\mathbb{P}'}_t|\leq 1+|\hat{a}^{1/2}_tz^{\mathbb{P}'}_t(t_2, Y_{t_2})|+|\hat{a}^{1/2}_tZ_t|.
$$
Defining 
$$
H^{\mathbb{P}'}_t:=\int^t_0\kappa^{\mathbb{P}'}_sdW^{\mathbb{P}'}_s,\ 0\leq t\leq t_2,
$$
we have
\begin{equation}\label{er2}||H^{\mathbb{P}'}||^2_{BMO(\mathbb{P}')}\leq C(1+||\hat{a}^{1/2}_tz^{\mathbb{P}'}_t(t_2, Y_{t_2})||^2_{\mathbb{H}^2_{BMO(\mathbb{P}')}}+||\hat{a}^{1/2}_tZ_t||^2_{\mathbb{H}^2_{BMO_2(\mathbb{P}')}}).\end{equation}
Applying a priori estimates for quadratic BSDEs (cf. Lemma 3.1 in Morlais \cite{M}), it is readily observed that
\begin{equation}\label{er}
||\hat{a}^{1/2}z^{\mathbb{P}'}(t_2, Y_{t_2})||^2_{H^2_{BMO(\mathbb{P}')}}
\leq C e^{4\gamma||Y||_{\mathbb{D}^\infty_H}}(1+||Y||_{\mathbb{D}^\infty_H}).
\end{equation}
Putting (\ref{er1}) and (\ref{er}) into (\ref{er2}), we deduce the following estimate uniformly in $\mathbb{P}'$:
\begin{equation}\label{es5}||H^{\mathbb{P}'}||^2_{BMO(\mathbb{P}')}\leq M_H,
\end{equation}
where $M_H$ depends simply on $||Y||_{\mathbb{D}^\infty_H}$. 
This implies that $H$ is a $BMO(\mathcal{P}_H)$-martingale.\\[6pt]
Define a probability measure $\mathbb{Q}'\ll\mathbb{P}'$ by
$\frac{d\mathbb{Q}'}{d\mathbb{P}'}|_{\mathcal{F}_t}=\mathcal{E}(\int^\cdot_0\kappa^{\mathbb{P}'}_sdW^{\mathbb{P}'}_s)_t
$
and a process
$M_t:= \exp(\int^t_{t_1}\lambda_sds),\ t_1\leq t\leq t_2.$
Applying It\^o's formula to $M\delta Y$ under $\mathbb{Q}'$, we have
\begin{align}\label{est2}
\delta Y^{\mathbb{P}'}_{t_1}=\mathbb{E}^{\mathbb{Q}'}_{t_1}\bigg[\int^{t_2}_{t_1}M_tdK^{\mathbb{P}'}_t\bigg]
&\leq\mathbb{E}^{\mathbb{Q}'}_{t_1}[\sup_{t_1\leq t\leq t_2}(M_t)(K^{\mathbb{P}'}_{t_2}-K^{\mathbb{P}'}_{t_1})]\\&\leq e^\mu \mathbb{E}^{\mathbb{P}'}_{t_1}\bigg[\frac{\mathcal{E}(H^{\mathbb{P}'})_{t_2}}{\mathcal{E}(H^{\mathbb{P}'})_{t_1}}(K^{\mathbb{P}'}_{t_2}-K^{\mathbb{P}'}_{t_1})\bigg]\notag, \ \mathbb{P}-a.s..
\end{align}
Thanks to (\ref{es5}) and Lemma \ref{BMOP}, we can find uniformly for all $\mathbb{P}'\in\mathcal{P}_H(t^+_1, \mathbb{P})$ two constants $q>1$ and $C_{RH}>0$ such that 
\begin{align*}
\delta Y^{\mathbb{P}'}_{t_1}&\leq C_{RH}^{{1}/{q}}e^\mu\mathbb{E}^{\mathbb{P}'}_{t_1}[(K^{\mathbb{P}'}_{t_2}-K^{\mathbb{P}'}_{t_1})^{p}]^{{1}/{p}}\\[6pt]
&\leq C^{{1}/{q}}_{RH}e^\mu(\mathbb{E}^{\mathbb{P}'}_{t_1}[(K^{\mathbb{P}'}_{t_2}-K^{\mathbb{P}'}_{t_1})^{2p-1}])^{1/2p}\mathbb{E}^{\mathbb{P}'}_{t_1}
[(K^{\mathbb{P}'}_{t_2}-K^{\mathbb{P}'}_{t_1})]^{{1}/{2p}}\\[3pt]
&\leq C^{{1}/{q}}_{RH}e^\mu(\sideset{}{^\mathbb{P}}\esssup_{\mathbb{P}'\in\mathcal{P}_H(t^+_1,\mathbb{P})}\mathbb{E}^{\mathbb{P}'}_{t_1}[(K^{\mathbb{P}'}_{t_2}-K^{\mathbb{P}'}_{t_1})^{2p-1}])^{1/2p}\mathbb{E}^{\mathbb{P}'}_{t_1}
[(K^{\mathbb{P}'}_{t_2}-K^{\mathbb{P}'}_{t_1})]^{{1}/{2p}}\\
&\leq  C^{{1}/{q}}_{RH}C_{2p-1}^{{1}/{2p}}e^\mu\mathbb{E}^{\mathbb{P}'}_{t_1}
[(K^{\mathbb{P}'}_{t_2}-K^{\mathbb{P}'}_{t_1})]^{{1}/{2p}},\notag\ \mathbb{P}-a.s.,
\end{align*}
where 
$1/p+1/q=1$. 
Since $C_{RH}$ and $C_{2p-1}$ are independent of $\mathbb{P}'$, we can take essential infimum over all $\mathbb{P}'\in\mathcal{P}_H(t^+_1, \mathbb{P})$ on the left-hand side of the above inequality and deduce by the minimum condition (\ref{min}) that
\begin{align*}
Y_{t_1}&-\sideset{}{^\mathbb{P}}\esssup_{\mathbb{P}'\in\mathcal{P}_H(t^+_1, \mathbb{P})}y^{\mathbb{P}'}_{t_1}(t_2, Y_{t_2})\\
&\leq C^{{1}/{q}}_{RH} C_{2p-1}^{{1}/{2p}}e^\mu\sideset{}{^\mathbb{P}}\essinf_{\mathbb{P}'\in\mathcal{P}_H(t^+_1, \mathbb{P})}[(K^{\mathbb{P}'}_{t_2}-K^{\mathbb{P}'}_{t_1})]^{{1}/{2p}}=0,\ \mathbb{P}-a.s..
\end{align*}
From (\ref{est2}), it is easily observed that $\delta Y^{\mathbb{P}'}_{t_1}\geq 0$, $\mathbb{Q}'$-a.s. and thus, $\mathbb{P}-a.s.$, for all $\mathbb{P}'\in\mathcal{P}_H(t^+_1, \mathbb{P})$, which directly yields the reverse inequality of (\ref{zbzd}). The proof of (\ref{er12}) is complete.\hfill{}$\square$\\[6pt]
\subsection{A priori estimates}
We now give some a priori estimates for quadratic 2BSDEs:
\begin{lemma}
Let (A1)-(A5) hold. Assume that $\xi\in L^\infty_H$ and that $(Y, Z)\in\mathbb{D}^\infty_H\times\mathbb{H}^2_H$ is a solution to 2BSDE (\ref{key}). Then, there exists a $C>0$ such that
\begin{equation}\label{poer}
||Y||_{\mathbb{D}^\infty_H}\leq C(1+||\xi||_{L^\infty_H})\ and\ ||Z||^2_{\mathbb{H}^2_{BMO(\mathcal{P}_H)}}\leq Ce^{4\gamma||\xi||_{L^\infty_H}}(1+||\xi||_{L^\infty_H}).\end{equation}
\end{lemma}
\noindent \textbf{Proof:} From (\ref{18}) and a priori estimates for quadratic BSDEs, we deduce the left-hand side of (\ref{poer}), whereas the right-hand side comes after (\ref{er1}). \hfill{}$\square$
\begin{lemma}\label{apri}
Let (A1)-(A5) hold. Assume that $\xi^i\in L^\infty_H$ and that $(Y^i, Z^i)\in\mathbb{D}^\infty_H\times\mathbb{H}^2_H$, $i=1, 2$, are two solution to 2BSDE (\ref{key}). Denote
\begin{align*}\delta\xi:=\xi^1-\xi^2,\ \delta Y:=Y^1-Y^2,\ 
\delta Z:=Z^1-Z^2,\\
\delta K^\mathbb{P}:=
(K^1)^{\mathbb{P}}-(K^2){^\mathbb{P}}\ and\ \Delta \delta Y_t= \delta Y_t-\delta Y_{t^-},
\end{align*}
then we have the following estimates
$$
||\delta Y||_{\mathbb{D}^\infty_H}\leq C||\delta\xi||_{L^\infty_H};$$
$$
\mathbb{E}^\mathbb{P}_\tau\bigg[\int^1_\tau|\hat{a}^{1/2}_t\delta Z_t|^2dt\bigg]
\leq C||\delta\xi||^2_{L^\infty_H}
\sum^2_{i=1}(1+e^{4\gamma||\xi^i||_{L^\infty_H}})(1+||\xi^i||_{L^\infty_H}).
$$
\end{lemma}
\noindent\textbf{Proof:} Similar to (\ref{est2}),
we can easily obtain the first inequality. For the second one, we apply It\^o's formula to $\delta Y^2$, then we have for a fixed $\mathbb{P}\in\mathcal{P}_H$ and a $\tau\in\mathcal{T}^1_0$, 
%
\begin{align*}
\delta Y^2_\tau+\int^1_\tau|\hat{a}^{1/2}_t\delta Z_t|^2dt
&\leq\delta\xi^2+2\int^1_\tau\delta Y_t(\hat{F}_t(Y^1_t, Z^1_t)-\hat{F}_t(Y^2_t, Z^2_t))dt\\
-&2\int^1_\tau\delta Y_t\delta Z_tdB_t+2\int^1_\tau \delta Y_{t^-}\delta dK^\mathbb{P}_t-\sum_{\tau< t\leq 1}(\Delta \delta Y_t)^2,\ \mathbb{P}-a.s..
\end{align*}
Taking expectation on the left-hand side and by (A3) and (\ref{18}), we deduce
\begin{align*}
\mathbb{E}^\mathbb{P}_\tau\bigg[\int^1_\tau|\hat{a}^{1/2}_t\delta Z_t|^2dt\bigg]
&\leq||\delta\xi||^2_{L^\infty_H}+2||\delta Y||_{\mathbb{D}^\infty_H}\bigg(2\alpha+\beta\sum^2_{i=1}||Y^i||_{\mathbb{D}^\infty_H}+\frac{\gamma}{2}\sum^2_{i=1}||Z^i||_{H^2_{BMO(\mathcal{P}_H)}}\bigg)\\
&+2||\delta Y||_{\mathbb{D}^\infty_H}(\mathbb{E}^\mathbb{P}_\tau[(K^1)^\mathbb{P}_1-(K^1)^\mathbb{P}_\tau]+\mathbb{E}^\mathbb{P}_\tau[(K^2)^\mathbb{P}_1-(K^2)^\mathbb{P}_\tau])\\
&\leq C||\delta\xi||_{L^\infty_H}
\sum^2_{i=1}(1+e^{4\gamma||\xi^i||_{L^\infty_H}})(1+||\xi^i||_{L^\infty_H}).
\end{align*}
We complete the proof.
\hfill{}$\square$\\[6pt]
By either Theorem \ref{rep} or Lemma \ref{apri}, we deduce immediately the uniqueness of $Y$.
We observe that 
$d\langle Y, B\rangle_t=Z_td\langle B\rangle_t$, $0\leq t\leq 1$, $\mathcal{P}_H$-q.s., which implies the uniqueness of $Z$.
\section{Existence of solutions to 2BSDEs}
\noindent In this section, we provide the existence result for the 2BSDE (\ref{key}) under (A1)-(A5) by a pathwise construction introduced in Soner et al. \cite{STZ1} and \cite{STZ2} with the technique so-called regular conditional probability distribution (r.p.c.d.), which can be find in Stroock and Varadhan \cite{SV}. 
\subsection{Regular conditional probability measures}
For the convenience of the reader, we recall some notations of r.p.c.d. in Soner et al. \cite{STZ1}.
\begin{itemize}[leftmargin=*]
\item For each $t\in[0,1]$, let $\Omega^t:=\{\tilde{\omega}\in\mathcal{C}([t, 1], \mathbb{R}^d), \tilde{\omega}(t)=0\}$ be the shifted space, $B^t$ the shifted canonical process, $\mathcal{F}^t$ the shifted filtration generated by $B^t$.\\[0pt]
\item For each $0\leq s\leq t\leq 1$ and $\omega\in\Omega^s$, we define the shifted path $\tilde{\omega}\in\Omega^t$ by
$$\tilde{\omega}_u:=\omega_u-\omega_t,\ u\in[t, 1].$$
\item For each $0\leq s\leq t\leq 1$, $\omega\in\Omega^s$ and $\tilde{\omega}\in\Omega^t$, we  define the concatenation path $\omega\otimes_t\tilde{\omega}\in\Omega^s$ by
$$(\omega\otimes_t\tilde{\omega})_u:=\omega_u\textnormal{\textbf{1}}_{[s, t)}(u)+(\omega_t+\tilde{\omega}_u)\textnormal{\textbf{1}}_{[t, 1]}(u),\ u\in [s, 1].$$
\item For each $0\leq s\leq t\leq 1$, $\omega\in\Omega^s$ and an $\mathcal{F}^s_1$-measurable random variable $\xi$ on $\Omega^s$,  we define the shifted $\mathcal{F}^t_1$-measurable random variable $\xi^{t, \omega}$ on $\Omega^t$ by
$$\xi^{t, \omega}(\tilde{\omega}):=\xi(\omega\otimes_t\tilde{\omega}), \tilde{\omega}\in\Omega^t.$$
\item For each $0\leq s\leq t\leq 1$, the shifted process $X^{t, \omega}$ of an $\mathcal{F}^s$-progressively measurable $X$ is $\mathcal{F}^t$-progressively measurable.\\[0pt]
\item For each $t\in[0,1]$ and $\omega\in\Omega$, we define our shifted generator by
$$\hat{F}^{t, \omega}_s(\tilde{\omega}, y, z):=F_s(\omega\otimes_t\tilde{\omega}, y, z, \hat{a}^t_s(\tilde{\omega})),\ (s, \tilde{\omega})\in[t, 1]\times\Omega^t.$$
\item For each $t\in[0,1]$, $\mathcal{P}^t_H$ denotes the collection of all those $\mathbb{P}\in \overline{\mathcal{P}}^t_S$ such that $$\underline{a}^\mathbb{P}\leq\hat{a}^t_s
\leq\overline{a}^\mathbb{P}\ {\rm and}\ \hat{a}^t_s\in D_{F_s},\ \lambda\times\mathbb{P}-a.e.,
$$
for some $\underline{a}^\mathbb{P}$, $\overline{a}^\mathbb{P}\in\mathbb{S}_d^{>0}$ and all $(y,z)\in\mathbb{R}\times\mathbb{R}^d$.\\
\item For $t\in[0,1]$ and $\mathbb{P}\in\mathcal{P}_H$, the r.p.c.d $\mathbb{P}^\omega_t$ of $\mathbb{P}$ induces naturally a probability measure $\mathbb{P}^{t, \omega}$ on $(\Omega^t, \mathcal{F}^t_1)$ which satisfies that for each bounded and $\mathcal{F}_1$-measurable random variable $\xi$,
$$\mathbb{E}^{\mathbb{P}^\omega_t}[\xi]=\mathbb{E}^{\mathbb{P}^{t, \omega}}[\xi^{t, \omega}].$$
\item By Lemma 4.1 in \cite{STZ1}, $\mathbb{P}^{t, \omega}$ is an elements in $\mathcal{P}^t_H$ and for each $t\in[0,1]$ and $\mathbb{P}\in\mathcal{P}_H$, it holds for $\mathbb{P}$-a.s. $\omega\in\Omega$,
\begin{align*}
F_s(\omega\otimes_t\tilde{\omega}, y, z, \hat{a}^t_s(\tilde{\omega}))=F_s(\omega\otimes_t\tilde{\omega}, y, z&, \hat{a}_s(\omega\otimes_t\tilde{\omega})),\ \lambda\times\mathbb{P}^{t, \omega}-a.e..
\end{align*}
\end{itemize}
\subsection{Existence result}For some fixed $\omega\in\Omega$ and $0\leq t_1\leq t_2\leq 1$, we consider a quadratic BSDE of the following type on the shifted space $\Omega^{t_1}$ under each $\mathbb{P}^{t_1}\in\mathcal{P}^{t_1}_H$:
\begin{align}\label{BSDER}
y^{\mathbb{P}^{t_1}, t_1, \omega}_s=\eta^{t_1, \omega}&+\int^{t_2}_s \hat{F}^{t_1, \omega}_u(y^{\mathbb{P}^{t_1}, t_1, \omega}_u, z^{\mathbb{P}^{t_1}, t_1, \omega}_u)du\\
&-\int^{t_2}_sz^{\mathbb{P}^{t_1}, t_1, \omega}_udB^{t_1}_u,\  t_1\leq s\leq t_2,\ \mathbb{P}^{t_1}-a.s.,\notag
\end{align}
where $\eta\in L^\infty_H$ is a $\mathcal{F}_{t_2}$-measurable random variable. It is well known that (\ref{BSDER}) admits a unique solution $
(y^{\mathbb{P}^{t_1}, t_1, \omega}(t_2, \eta), z^{\mathbb{P}^{t_1}, t_1, \omega}(t_2, \eta))$ under (A1)-(A5).
In view of the Blumenthal zero-one law, $y^{\mathbb{P}^{t_1}, t_1, \omega}_{t_1}(t_2, \eta)$ is a deterministic constant $\mathbb{P}^{t_1}$-a.s. for any given $\eta$ and $\mathbb{P}^{t_1}$. \\[6pt]
The following lemma describes the relationship between $y_t^\mathbb{P}(1, \xi)$ and $y^{\mathbb{P}^{t, \omega}, t, \omega}_t(1, \xi)$, where the first one is the solution to (\ref{st}) with the parameters $(1, \xi)$ under a fixed $\mathbb{P}\in\mathcal{P}_H$
and the second one is the solution to (\ref{BSDER}) when $t_1$ is a fixed $t$, $\mathbb{P}^{t}$ is in fact the r.p.c.d $\mathbb{P}^{t, \omega}$ of $\mathbb{P}$ and $(t_2, \eta)=(1, \xi)$.
\begin{lemma}\label{keylemma}
Assume (A1)-(A5) hold. For a given $\xi\in L^\infty_H$ and a fixed $\mathbb{P}\in\mathcal{P}_H$, 
we have, for each $t\in[0,1]$ and $\mathbb{P}$-a.s. $\omega\in\Omega$,
\begin{equation}\label{equa}
y^\mathbb{P}_t(1, \xi)(\omega)=y^{\mathbb{P}^{t, \omega}, t, \omega}_t(1, \xi).\end{equation}
\end{lemma}
\noindent\textbf{Proof:} 
Similar to the proof of Lemma 4.1 in Soner et al. \cite{STZ1}, we have the following conclusion: because $\xi\in L^\infty(\mathbb{P})$, for $\mathbb{P}$-a.s. $\omega\in\Omega$, $|\xi^{t, \omega}|\leq ||\xi||_{L^\infty(\mathbb{P})}$, $\mathbb{P}^{t, \omega}$-a.s..
Thus, (\ref{BSDER}) is well defined under our setting and the right-hand side of (\ref{equa}) is the unique solution of (\ref{BSDER}).\\[6pt]
\noindent We emphasize that the wellposedness of both (\ref{st}) and (\ref{BSDER}) as well as the estimates of the solutions are already provided by Kobylanski \cite {KO} and Morlais \cite{M}. Our job here is only to redo the construction of two sequences formed by the solutions of Lipschitz BSDEs, which approximate the solutions on both sides of (\ref{equa}).\\[6pt]
\noindent By Lemma3.1 in Morlais \cite{M}, we can find an $M:=e^\beta(\alpha+
||\xi(\omega)||_{L^\infty(\mathbb{P})}),$ which is the bound of both side of (\ref{equa}). Then, we choose a $\mathcal{C}^{\infty}(\mathbb{R})$ function which takes value in $[0, 1]$ and satisfies that
$$
\phi(u)=\left\{
\begin{array}{l@{\quad, \quad}l}
1 & u\in[e^{-\gamma M}, e^{\gamma M}];\\
0 & u\in(-\infty, e^{-\gamma (M+1)}]\cup [e^{\gamma (M+1)}, +\infty).
\end{array}\right.\\
$$
We can verify that for each $t\in[0, 1]$, 
\begin{equation}\label{wuliao1}
\mathcal{Y}^\mathbb{P}_t(1, e^{\gamma\xi}, \hat{G}):=\exp(\gamma y^\mathbb{P}_t(1, \xi)),
\end{equation}
solves a quadratic BSDE with the parameters $(1, e^{\gamma\xi})$ and the generator $\hat{G}$ of the following form: for each $(\omega, t, \mathcal{Y}, \mathcal{Z})\in\Omega\times[0,1]\times\mathbb{R}\times\mathbb{R}^d$,
\begin{equation}\label{geng}\hat{G}_t(\omega, \mathcal{Y}, \mathcal{Z})
:=\phi(\mathcal{Y})\bigg(\gamma\mathcal{Y}
\hat{F}_t\bigg(\omega, \frac{\ln(\mathcal{Y})}{\gamma}, \frac{\mathcal{Z}}{\gamma\mathcal{Y}}\bigg)
-\frac{1}{2\mathcal{Y}}|\hat{a}^{1/2}_t(\omega)\mathcal{Z}|^2\bigg).
\end{equation}
On the other hand, fixing $(\omega, t)\in\Omega\times[0,1]$,
\begin{equation}\label{wuliao3}
\mathcal{Y}^{\mathbb{P}^{t,\omega}, t, \omega}_s(1, e^{\gamma\xi}, \hat{G}^{t, \omega})(\tilde{\omega}):=\exp(\gamma y^{\mathbb{P}^{t,\omega}, t, \omega}_s(1, \xi)(\tilde{\omega})),\ t\leq s\leq 1,
\end{equation}
defines a solution that solves a quadratic BSDE under $\mathbb{P}^{t, \omega}$ with the parameters $(1, e^{\gamma\xi})$ and the generator $\hat{G}^{t, \omega}$ of the following form:  
for each $(\tilde{\omega}, s, \mathcal{Y}, \mathcal{Z})\in\Omega^t\times[t,1]\times\mathbb{R}\times\mathbb{R}^d$,
$$\hat{G}^{t, \omega}_s(\tilde{\omega}, \mathcal{Y}, \mathcal{Z})
:=\phi(\mathcal{Y})\bigg(\gamma\mathcal{Y}
\hat{F}_s^{t, \omega}\bigg(\tilde{\omega}, \frac{\ln(\mathcal{Y})}{\gamma}, \frac{\mathcal{Z}}{\gamma\mathcal{Y}}\bigg)
-\frac{1}{2\mathcal{Y}}|(\hat{a}^t_s)^{1/2}(\tilde{\omega})\mathcal{Z}|^2\bigg).
$$
Now, our main aim is changed into that for each $t\in[0,1]$ and $\mathbb{P}$-a.s. $\omega\in\Omega$,
\begin{equation}\notag
\mathcal{Y}^\mathbb{P}_t(1, e^{\gamma\xi}, \hat{G})(\omega)=\mathcal{Y}^{\mathbb{P}^{t, \omega}, t, \omega}_t(1, e^{\gamma\xi}, \hat{G}^{t, \omega}).
\end{equation}
\noindent For each $(\omega, t, \mathcal{Y}, \mathcal{Z})\in\Omega\times[0,1]\times\mathbb{R}\times\mathbb{R}^d$, we set
\begin{equation}\label{hatg}
\hat{G}^n_t(\omega, \mathcal{Y}, \mathcal{Z}):=\sup_{(p, q)\in\mathbb{Q}\times\mathbb{Q}^d}\{\hat{G}_t(\omega, p, q)-n|(p-\mathcal{Y})|-n|\hat{a}^{1/2}_t(\omega)(q-\mathcal{Z})|\},\ n\in\mathbb{N}
\end{equation}
and also for fixed $(\omega, t)\in\Omega\times[0,1]$ and each $(\tilde{\omega}, s, \mathcal{Y}, \mathcal{Z})\in\Omega^t\times[t,1]\times\mathbb{R}\times\mathbb{R}^d$, we define
\begin{align*}((\hat{G}^{t, \omega})^n_s&(\tilde{\omega}, \mathcal{Y}, \mathcal{Z})\\
&:=\sup_{(p, q)\in\mathbb{Q}\times\mathbb{Q}^d}\{\hat{G}^{t, \omega}_s(\tilde{\omega}, p, q)-n|(p-\mathcal{Y})|-n|(\hat{a}^t_s)^{1/2}(\tilde{\omega})(q-\mathcal{Z})|\},\ n\in\mathbb{N}.
\end{align*}
By Lemma 4.1 in Soner et al. \cite{STZ1}, for each $t\in[0,1]$, $\mathbb{P}$-a.s. $\omega\in\Omega$ and each $(y, z)\in\mathbb{R}\times\mathbb{R}^d$, 
\begin{align}\label{global}
\hat{F}^{t, \omega}_s(\tilde{\omega}, y, z)&=F^{t, \omega}_s(\tilde{\omega}, y, z, \hat{a}^t_s(\tilde{\omega}))=F_s(\omega\otimes_t\tilde{\omega}, y, z, \hat{a}_s(\omega\otimes_t\tilde{\omega}))\\
&=(\hat{F}(\cdot, \cdot))^{t, \omega}_s(\tilde{\omega}, y, z),
\ {\rm and}\ \hat{a}^t_s(\tilde{\omega})=\hat{a}^{t, \omega}_s(\tilde{\omega}),
\ \lambda\times\mathbb{P}^{t, \omega}-a.e..\notag
\end{align}
We call $(\hat{F}(\cdot, \cdot))^{t, \omega}$ the globally shifted generator of $\hat{F}$.
From (\ref{global}), we can deduce that for each $t\in[0,1]$, $\mathbb{P}$-a.s. $\omega\in\Omega$ and each $(\mathcal{Y}, \mathcal{Z})\in\mathbb{R}\times\mathbb{R}^d$, 
\begin{equation*}\label{wuliao2}
\hat{G}^{t, \omega}_s(\tilde{\omega}, \mathcal{Y}, \mathcal{Z})=(\hat{G}(\cdot, \cdot))^{t, \omega}_s(\tilde{\omega}, \mathcal{Y}, \mathcal{Z}),
\ \lambda\times\mathbb{P}^{t, \omega}-a.e.,
\end{equation*}
and furthermore that for each $n\in\mathbb{N}$,
\begin{equation}\label{keyrel}
(\hat{G}^{t, \omega})^n_s(\tilde{\omega}, \mathcal{Y}, \mathcal{Z})
=(\hat{G}^n(\cdot, \cdot))^{t, \omega}_s(\tilde{\omega}, \mathcal{Y}, \mathcal{Z}),\ \lambda\times\mathbb{P}^{t, \omega}-a.e..
\end{equation}
Moreover, it is easy to verify that  for each $(\omega, t, \mathcal{Y}, \mathcal{Z})\in\Omega\times[0,1]\times\mathbb{R}\times\mathbb{R}^d$,
\begin{align*}-e^{M+1}(\alpha\gamma+\beta(M+1))&-e^{M+1}|\hat{a}_t^{1/2}(\tilde{\omega})\mathcal{Z}|^2\leq
\hat{G}_t(\omega, \mathcal{Y}, \mathcal{Z})\\ &\leq \hat{G}^{n+1}_t(\omega, \mathcal{Y}, \mathcal{Z})\leq\hat{G}^n_t(\omega, \mathcal{Y}, \mathcal{Z})\leq e^{M+1}(\alpha\gamma+\beta(M+1)),
\end{align*}
and 
$\hat{G}^{n}_t(\omega, \mathcal{Y}, \mathcal{Z})\downarrow \hat{G}_t(\omega, \mathcal{Y}, \mathcal{Z})
$
uniformly on compact sets in $[0, 1]\times\mathbb{R}\times\mathbb{R}^d$. Similarly,  for fixed $(\omega, t)\in\Omega\times[0,1]$ and each $(\tilde{\omega}, s, \mathcal{Y}, \mathcal{Z})\in\Omega^t\times[t,1]\times\mathbb{R}\times\mathbb{R}^d$, 
\begin{align*}-e^{M+1}(\alpha\gamma+\beta(M+1))&-e^{M+1}|(\hat{a}^t_s)^{1/2}(\tilde{\omega})\mathcal{Z}|^2\leq
\hat{G}_t(\omega, \mathcal{Y}, \mathcal{Z})\\ &\leq (\hat{G}^{t, \omega})^{n+1}_t(\omega, \mathcal{Y}, \mathcal{Z})\leq(\hat{G}^{t, \omega})^n_t(\omega, \mathcal{Y}, \mathcal{Z})\leq e^{M+1}(\alpha\gamma+\beta(M+1)),
\end{align*}
and 
$(\hat{G}^{t, \omega})^{n}_t(\omega, \mathcal{Y}, \mathcal{Z})\downarrow (\hat{G}^{t, \omega})_t(\omega, \mathcal{Y}, \mathcal{Z})
$
uniformly on compact sets in $[t, 1]\times\mathbb{R}\times\mathbb{R}^d$. By Lemma 3.3 (monotone stability) in Morlais \cite{M}, we have, for each $t\in[0,1]$ and $\mathbb{P}$-a.s. $\omega\in\Omega$,  
\begin{equation}\label{o2}
\mathcal{Y}_t^\mathbb{P}(1, e^{\gamma\xi}, \hat{G}^n)(\omega)\downarrow\mathcal{Y}_t^\mathbb{P}(1, e^{\gamma\xi}, \hat{G})(\omega),\ {\rm as}\ n\rightarrow+\infty,
\end{equation}
and for fixed $(\omega, t)\in\Omega\times [0, 1]$, 
\begin{equation}\label{o1}
\mathcal{Y}_t^{\mathbb{P}^{t, \omega}, t, \omega}(1, e^{\gamma\xi}, (\hat{G}^{t, \omega})^n)\downarrow\mathcal{Y}_t^{\mathbb{P}^{t, \omega}, t, \omega}(1, e^{\gamma\xi}, \hat{G}^{t, \omega}),\ {\rm as}\
n\rightarrow+\infty.
\end{equation}
To obtain the desired result, it is suffice to prove that for each $n\in\mathbb{N}$, a fixed $t\in[0,1]$ and $\mathbb{P}$-a.s. $\omega\in\Omega$, 
\begin{equation}\label{FEFE}
\mathcal{Y}^\mathbb{P}_t(1, e^{\gamma\xi}, \hat{G}^n)(\omega)=\mathcal{Y}^{\mathbb{P}^{t, \omega}, t, \omega}_t(1, e^{\gamma\xi}, (\hat{G}^{t, \omega})^n).
\end{equation}
We notice that the generators of both sides of (\ref{FEFE}) satisfy the following uniform Lipschitz conditions: for each $n\in\mathbb{N}$ and $(\omega, t, \mathcal{Y}^1, \mathcal{Y}^2, \mathcal{Z}^1, \mathcal{Z}^2)\in\Omega\times[0,1]\times\mathbb{R}\times\mathbb{R}\times\mathbb{R}^d\times\mathbb{R}^d$,
$$
|\hat{G}^n_t(\omega, \mathcal{Y}^1, \mathcal{Z}^1)-\hat{G}^n_t(\omega, \mathcal{Y}^2, \mathcal{Z}^2)|\leq n|\mathcal{Y}^1-\mathcal{Y}^2|+n|\hat{a}^{1/2}_t(\omega)(\mathcal{Z}^1-\mathcal{Z}^2)|;
$$
and for fixed $(\omega, t)\in\Omega\times[0,1]$, each $n\in\mathbb{N}$ and
$(\tilde{\omega}, s, \mathcal{Y}^1, \mathcal{Y}^2, \mathcal{Z}^1, \mathcal{Z}^2)\in\Omega^t\times[t,1]\times\mathbb{R}\times\mathbb{R}\times\mathbb{R}^d\times\mathbb{R}^d$,
$$
|(\hat{G}^{t,\omega})^n_s(\tilde{\omega}, \mathcal{Y}^1, \mathcal{Z}^1)-(\hat{G}^{t, \omega})^n_s(\tilde{\omega}, \mathcal{Y}^2, \mathcal{Z}^2)|\leq n|\mathcal{Y}^1-\mathcal{Y}^2|+n|(\hat{a}^t_s)^{1/2}(\tilde{\omega})(\mathcal{Z}^1-\mathcal{Z}^2)|.
$$
Since the solutions of these Lipschitz BSDEs can be constructed via Picard iteration, from (\ref{keyrel}), we can obtain (\ref{FEFE})
(cf. (i) in the proof of Proposition 4.7 in Soner et al. \cite{STZ1} and (i) in the proof of Proposition 5.1 in Possamai and Zhou \cite{PZ}). 
Then, (\ref{o2}) and (\ref{o1}) give the desired result.\hfill{}$\square$
\begin{remark}
The lemma above is the key point of this paper, which ensure us to prove the following proposition under the Kobylanski \cite{KO} type condition (A3) instead of Tevzadze \cite{TEV} type one, which is adopted by Possamai and Zhou \cite{PZ} to make sure that the solutions of quadratic BSDEs on both original and shifted spaces can be constructed via Picard iteration, so that the statement corresponding to (\ref{equa}) in Soner et al. \cite{STZ1} for Lipschitz BSDEs still holds. In the present paper, the lemma above is proved by a monotonic convergence technique for classical BSDEs under a fixed $\mathbb{P}$, but it is still difficult to obtain a globally monotonic convergence theorem for quadratic 2BSDEs, 
as Possamai and Zhou \cite{PZ} has already stated.
\end{remark}
\noindent Similar to the one in Soner et al. \cite{STZ1}, we define the following value process $V_t$ pathwisely: for each $(\omega, t)\in\Omega\times[0,1]$,
\begin{equation}\label{v}
V_t(\omega):=\sup_{\mathbb{P}^t\in\mathcal{P}^t_H}y^{\mathbb{P}^t, t, \omega}_t(1, \xi).
\end{equation}
For the rest part of the proof of the existence, we assume moreover that
the terminal value $\xi$ is an element in $UC_b(\Omega)$. Therefore,
it is readily observed that for all $(\omega, t)\in\Omega\times[0,1]$,
\begin{equation}\label{estv}
V_t(\omega)\leq C(1+\sup_{\omega\in\Omega}|\xi(\omega)|),
\end{equation}
and 
there exists a modulus of continuity $\rho$, such that for each $t\in[0, 1]$ and $(\omega, \omega', \tilde{\omega})\in \Omega\times\Omega\times\Omega^t
$, 
$$
|\xi^{t, \omega}(\tilde{\omega})-\xi^{t, \omega'}(\tilde{\omega})|\leq \rho(||\omega-\omega'||^\infty_t).
$$
\noindent Recalling the uniform continuity of $F$ in $\omega$,
we have moreover that for each $0\leq t\leq s\leq 1$, $(\omega, \omega', \tilde{\omega}, y, z)\in\Omega\times\Omega\times\Omega^t\times\mathbb{R}\times\mathbb{R}^d$,
$$
|\hat{F}^{t, \omega}_s(\tilde{\omega}, y, z)-\hat{F}^{t, \omega'}_s(\tilde{\omega}, y, z)|\leq \rho(||\omega-\omega'||^\infty_t).
$$
We define
$$\delta y:=y^{\mathbb{P}^t, t, \omega}(1, \xi)-y^{\mathbb{P}^t, t, \omega'}(1, \xi),\ \delta z:=z^{\mathbb{P}^t, t, \omega}(1, \xi)-z^{\mathbb{P}^t, t, \omega'}(1, \xi), $$
$$
\delta\xi:=\xi^{t, \omega}-\xi^{t, \omega'},\ \delta\hat{F}(y, z):=\hat{F}^{t, \omega}(y, z)-\hat{F}^{t, \omega'} (y, z).$$
Proceeding the same in the proof of Theorem \ref{rep}, for each $(\omega, \omega', t)\in\Omega\times\Omega\times [0,1]$ and a fixed $\mathbb{P}^t\in\mathcal{P}^t_H$, we can find a $\mathbb{Q}^t\ll\mathbb{P}^t$ and a bounded process $M$, such that
\begin{equation}\label{deva}
|\delta y_t|=\mathbb{E}^{\mathbb{Q}^t}\bigg[M_1\delta\xi
+\int^1_tM_s\delta\hat{F}_s(y_s^{\mathbb{P}^t, t, \omega}(1, \xi), z_s^{\mathbb{P}^t, t, \omega}(1, \xi))ds\bigg]
\leq C\rho(||\omega-\omega'||^\infty_t).
\end{equation}
By the arbitrariness of $\mathbb{P}^t$, it follows that
\begin{equation}\label{dev}|V_t(\omega)-V_t(\omega')|\leq C\rho(||\omega-\omega'||^\infty_t),
\end{equation}
from which we can deduce that $V_t\in\mathcal{F}_t$.\\[6pt]
Parallel to Proposition 4.7  in Soner et al. \cite{STZ1} and Proposition 5.1 in Possamai and Zhou \cite{PZ}, we give the following dynamic programming principle:
\begin{proposition}\label{repv}
Under (A1)-(A5) and for a given $\xi\in UC_b$, we have, for each
$0\leq t_1\leq t_2\leq 1$ and $\omega\in\Omega$,
\begin{equation}\label{TH}
V_{t_1}(\omega)=\sup_{\mathbb{P}^{t_1}\in\mathcal{P}^{t_1}_H}y^{\mathbb{P}^{t_1}, t_1, \omega}_{t_1}(t_2, V_{t_2}).
\end{equation}
\end{proposition}
\noindent\textbf{Proof:} Without loss of generality, we only need to prove the case when $t_1=0$ and $t_2=t$, i.e.,
$$V_0=\sup_{\mathbb{P}\in\mathcal{P}_H}y^\mathbb{P}_0(t, V_t).$$
Fixing $\mathbb{P}\in\mathcal{P}_H$, for each $\omega\in\Omega$ and $t\in[0, 1]$, $\mathbb{P}^{t, \omega}\in\mathcal{P}^t_H$. By Lemma \ref{keylemma} and from (\ref{v}), we have, for each $t\in [0,1]$ and $\mathbb{P}$-a.s. $\omega\in\Omega$,
\begin{equation*}
y^\mathbb{P}_t(1, \xi)(\omega)=y^{\mathbb{P}^{t, \omega}, t, \omega}_t(1, \xi)\leq \sup_{\mathbb{P}^t\in\mathcal{P}^t_H}y^{\mathbb{P}^t, t, \omega}_t(1, \xi)=V_t(\omega).
\end{equation*}
Applying Theorem 2.7 (comparison principle) in Morlais \cite{M},
it follows that $V_0\leq\sup_{\mathbb{P}\in\mathcal{P}_H}y^\mathbb{P}_0(t, V_t)$. \\[6pt]
We omit the rest part of the proof, i.e., the proof of the reverse inequality of (\ref{TH}) via the r.c.p.d. technique, since it goes in a exact same way as the one in Soner et al. \cite{STZ1} and Possamai and Zhou \cite{PZ}.
\hfill{}$\square$\\[6pt]
We shall head to the Doob-Meyer type decomposition of $V$ based on some results for quadratic $g$-supermartingales in Ma and Yao \cite{MY}. These results were obtained under the assumptions for the proof of uniqueness in Kobylanski \cite{KO}, since these assumptions ensure that the wellposedness of the corresponding quadratic BSDEs, so that the $g$-expectation can be well defined. 
However, Morlais \cite{M} also provided the wellposedness of quadratic BSDEs under (A3)-(A5) for the BSDEs of the form (\ref{st}), then the applicability of these arguments to such type of $\hat{F}$-supermartingales will not alter under each $\mathbb{P}\in\mathcal{P}_H$.\\[6pt]
%
For a fixed $\mathbb{P}\in\mathcal{P}_H$, from (\ref{TH}) and by Lemma \ref{keylemma}, we have for each $0\leq t_1\leq t_2\leq 1$,
\begin{equation}\label{fefe}
V_{t_1}\geq y^\mathbb{P}_{t_1}(t_2, V_{t_2}),\ \mathbb{P}-a.s..
\end{equation}
Thus, by Definition 5.1 in Ma and Yao \cite{MY}, $V$ is an $\hat{F}$-supermartingale under $\mathbb{P}$. Then, for each $(\omega, t)\in\Omega\times[0,1]$, we define 
\begin{equation*}
V^+_t(\omega):=\limsup_{\mathbb{Q}\cap(t, 1]\ni r\downarrow t} V_r(\omega).
\end{equation*}
Applying corollary 5.6 (downcrossing inequality) in Ma and Yao \cite{MY} , one can see that for $\mathbb{P}$-a.s. $\omega\in\Omega$, $\lim_{\mathbb{Q}\cap(t, 1]\ni r\downarrow t} V_r$ exists for all $t\in[0, 1]$. Therefore, we have
\begin{equation}\label{v+}
V^+_t=\lim_{\mathbb{Q}\cap(t, 1]\ni r\downarrow t} V_r,\ 0\leq t\leq 1,\ \mathcal{P}_H-q.s.,
\end{equation}
which implies that $V^+$ has $\mathcal{P}_H$-q.s. c\`adl\`ag paths.\\[6pt]
The following proposition (corresponding to Proposition 4.10 and 4.11 in Soner et al. \cite {STZ1}  and Proposition 5.2 in Possamai and Zhou \cite{PZ}) demonstrates the relationship between $V$ and $V^+$, from the second part of which, we can deduce that $V$ is a c\`adl\`ag $\hat{F}$-supermartingale under each $\mathbb{P}\in\mathcal{P}_H$, then we apply the Doob-Meyer type decomposition theorem (cf. Theorem 5.8 in Ma and Yao \cite{MY}) directly to $V$.
\begin{proposition}\label{vv}
Assume (A1)-(A5) hold. For a given $\xi\in UC_b(\Omega)$ and a fixed $\mathbb{P}\in\mathcal{P}_H$, we define
\begin{equation}\label{vp+}
V^\mathbb{P}_t:=\sideset{}{^\mathbb{P}}\esssup_{\mathbb{P}'\in\mathcal{P}_H(t, \mathbb{P})}y^{\mathbb{P}'}_t(1, \xi)\ {and}\ V^{\mathbb{P}, +}_t:=\sideset{}{^\mathbb{P}}\esssup_{\mathbb{P}'\in\mathcal{P}_H(t^+, \mathbb{P})}y^{\mathbb{P}'}_t(1, \xi).
\end{equation}
Then, we have
$$V_t=V^{\mathbb{P}}_t\ { and}\ V^+_t=V^{\mathbb{P}, +}_t,\ 0\leq t\leq 1,\ \mathbb{P}-a.s..$$
Moreover, 
\begin{equation}\label{vp++}V_t=V^+_t,\ 0\leq t\leq 1,\ \mathcal{P}_H-q.s..\end{equation}
\end{proposition}
\noindent\textbf{Proof:}
For the proof of the first equality in (\ref{vp+}) and that $V^+_t\geq V^{\mathbb{P}, +}_t$, we can proceed the same steps in the proof of Proposition 4.10 in Soner et al. \cite{STZ1}.  Here, we would like only to prove that for a fixed $\mathbb{P}\in\mathcal{P}_H$, 
$$V^+_t\leq V^{\mathbb{P},+}_t,\ 0\leq t\leq 1,\ \mathbb{P}-a.s.,$$
since the technique will be a little different. \\[6pt]
Fixing a $\mathbb{P}\in\mathcal{P}_H$ , a $t\in[0, 1]$ and an $r\in\mathbb{Q}\cap(t, 1]$, from the first equality, we have
$$
V^\mathbb{P}_r:=\sideset{}{^\mathbb{P}}\esssup_{\mathbb{P}'\in\mathcal{P}_H(t, \mathbb{P})}y^{\mathbb{P}'}_r(1, \xi).
$$
Following step 3 in the proof of Theorem 4.3 in Soner et al. \cite{STZ2}, we could find a sequence of probability measures such that  $\{\mathbb{P}_n\}_{n\in\mathbb{N}}\subset\mathcal{P}_H(r, \mathbb{P})\subset\mathcal{P}_H(t^+, \mathbb{P})$ and $y^{\mathbb{P}_n}_r(1, \xi)\uparrow V_r$, $\mathbb{P}$-a.s.. We consider the following BSDE with the parameters $(r, V_r)$ and the generator $\hat{F}$:
$$
y^\mathbb{P}_s=V_r+\int^r_s\hat{F}_u(y^\mathbb{P}_u, z^\mathbb{P}_u)du-\int^r_sz^\mathbb{P}_udB_u,\ 0\leq s\leq r,\ \mathbb{P}-a.s..
$$
and denote by $(y^\mathbb{P}(r, V_r), z^\mathbb{P}(r, V_r)))$ its solution. Then, it follows by Lemma 3.3 (monotone stability) in Morlais \cite{M} that
$$
y^\mathbb{P}_t(r, V_r)=y^\mathbb{P}_t(r, \lim_{n\rightarrow+\infty}y^{\mathbb{P}_n}_r(1, \xi))=
\lim_{n\rightarrow+\infty}y^{\mathbb{P}^n}_t(1, \xi)\leq V^{\mathbb{P}, +}_t,\ \mathbb{P}-a.s..
$$
Now, our aim is to find a sequence $\{r_m\}_{m\in\mathbb{N}}\subset (t, 1]$ such that $r_m\downarrow t$ and
\begin{equation}\label{sta34}
\lim_{m\rightarrow +\infty}y^\mathbb{P}_t(r_m, V_{r_m})= V^+_t,\ \mathbb{P}-a.s..
\end{equation}
Noticing that 
the generator $\hat{F}$ is no longer Lipschitz in $|\hat{a}^{1/2}Z|$, in general, the statement above (\ref{sta34}) is not straightforward if only the conditions that $V$ is uniformly bounded on $(t, 1]$ and that $V_r\rightarrow V^+_t$, $\mathbb{P}$-a.s. are given.
We define a sequence of BSDEs under $\mathbb{P}$ with the parameters $(r, e^{\gamma V_r})$ and the generator $
\hat{G}$ in the form of (\ref{geng}) and denote by $(\mathcal{Y}^\mathbb{P}(r, e^{\gamma V_r}, \hat{G}), \mathcal{Z}^\mathbb{P}(r, e^{\gamma V_r},\hat{G}))$ their solution. From the relationship that 
$$
\mathcal{Y}^\mathbb{P}_t(r, e^{\gamma V_r}, \hat{G})=e^{\gamma y^\mathbb{P}_t(r, V_r)},\ 0\leq t\leq r,\ \mathbb{P}-a.s.,
$$
the statement (\ref{sta34}) is equivalent to
\begin{equation}\label{sta35}
\lim_{m\rightarrow +\infty}\mathcal{Y}^\mathbb{P}_t(r_m, e^{\gamma V_{r_m}}, \hat{G})= e^{\gamma V^+_{t}},\ \mathbb{P}-a.s..
\end{equation}
Proceeding the same in the proof of Lemma \ref{keylemma}, for each $n$, we consider the solutions $\mathcal{Y}^{\mathbb{P}}(r, e^{\gamma V_r}, \hat{G}^n)$ of the BSDE with parameters $(r, e^{\gamma V_r})$ and the generator $\hat{G}^n$ in the form of (\ref{hatg}). Note $M:=C(1+\sup_{\omega\in\Omega}e^{\xi(\omega)})$ that is the uniform bound for all $\mathcal{Y}^{\mathbb{P}}(r, e^{\gamma V_r}, \hat{G}^n)$, we have 
$$
\mathbb{E}^\mathbb{P}[|\mathcal{Y}^{\mathbb{P}}(r, e^{\gamma V_r}, \hat{G}^n)-e^{\gamma V_r}|^2]\leq 
C(1+n+\alpha_M)(r-t),
$$
where $n$ is the Lipschitz constant of $\hat{G}^n$ and $\alpha_M>0$ depends simply on $M$. Therefore, for a fixed $n\in\mathbb{N}$, there exists a sequence $\{r^n_m\}_{m\in\mathbb{N}}\subset (t, 1]$ such that $r^n_m\downarrow t$ and 
\begin{equation*}\label{sta35}
\lim_{m\rightarrow +\infty}|\mathcal{Y}^{\mathbb{P}}_t(r^n_m, e^{\gamma V_{r^n_m}}, \hat{G}^n)-e^{\gamma V_{r_m}}|=0,\ \mathbb{P}-a.s.,
\end{equation*}
which implies
\begin{align*}
\lim_{m\rightarrow +\infty}&|\mathcal{Y}^{\mathbb{P}}_t(r^n_m, e^{\gamma V_{r^n_m}}, \hat{G}^n)-e^{\gamma V^+_{t}}|\\&\leq \lim_{m\rightarrow +\infty}|\mathcal{Y}^{\mathbb{P}}_t(r^n_m, e^{\gamma V_{r^n_m}}, \hat{G}^n)-e^{\gamma V_{r^n_m}}|+\lim_{m\rightarrow +\infty}|e^{\gamma V_{r^n_m}}-e^{\gamma V^+_{t}}|=0,\ \mathbb{P}-a.s..
\end{align*}
By the diagonal argument, we could find a universal sequence $\{\tilde{r}_m\}_{m\in\mathbb{N}}\subset (t, 1]$ such that $\tilde{r}_m\downarrow t$ and for each $n\in\mathbb{N}$, 
$$
\lim_{m\rightarrow +\infty}|\mathcal{Y}^{\mathbb{P}}_t(\tilde{r}_m, e^{\gamma V_{\tilde{r}_m}}, \hat{G}^n)-e^{\gamma V^+_{t}}|=0,\ \mathbb{P}-a.s..
$$
For each $n$, $m\in\mathbb{N}$,
$$
|\mathcal{Y}^{\mathbb{P}}_t(\tilde{r}_m, e^{\gamma V_{\tilde{r}_m}}, \hat{G}^n)|\leq M,
$$
and Lemma 3.3 (monotone stability) in Morlais \cite{M} shows that for each
$m\in\mathbb{N}$,
the following statement holds true $\mathbb{P}$-a.s.:
$$
\mathcal{Y}^{\mathbb{P}, n}_t(\tilde{r}_m, e^{\gamma V_{\tilde{r}_m}}, \hat{G}^n)\downarrow \mathcal{Y}^{\mathbb{P}}_t(\tilde{r}_m, e^{\gamma V_{\tilde{r}_m}}), \ as\ n\rightarrow +\infty.
$$
Thus,
\begin{align*}
\lim_{m\rightarrow +\infty}\mathcal{Y}^{\mathbb{P}}_t(\tilde{r}_m, e^{\gamma V_{\tilde{r}_m}})
&=\lim_{m\rightarrow +\infty}\lim_{n\rightarrow +\infty}\mathcal{Y}^{\mathbb{P}}_t(\tilde{r}_m, e^{\gamma V_{\tilde{r}_m}}, \hat{G}^n)\\&=\lim_{n\rightarrow +\infty}\lim_{m\rightarrow +\infty}\mathcal{Y}^{\mathbb{P}}_t(\tilde{r}_m, e^{\gamma V_{\tilde{r}_m}}, \hat{G}^n)=\lim_{n\rightarrow +\infty}e^{\gamma V^+_t}=e^{\gamma V^+_t},
\end{align*}
which ends the proof of (\ref{vp+}). Subsequently, the statement (\ref{vp++}) could be proved in a similar way as Proposition 4.11 in Soner et al. \cite{STZ1}. \hfill{}$\square$
\begin{theorem}
Under (A1)-(A5) and for a given $\xi\in UC_b(\Omega)$, the 2BSDE (\ref{key}) has a unique solution $(Y, Z)\in \mathbb{D}^{\infty}_H\times\mathbb{H}^2_H$.
\end{theorem}
\noindent\textbf{Proof:} From (\ref{fefe}) and (\ref{vp++}), we know that $V$ is a c\`adl\`ag $\hat{F}$-supermartingale. Applying the Doob-Meyer type decomposition (cf. Theorem 5.8 in Ma and Yao \cite{MY}) under each $\mathbb{P}\in\mathcal{P}_H$,
\begin{equation*}
V_t=V_1+\int^1_t \hat{F}_s(V_s, Z^\mathbb{P}_s)ds-\int^1_tZ^\mathbb{P}_sdBs+K^\mathbb{P}_1-K^\mathbb{P}_t,\ 0\leq t\leq1,\ \mathbb{P}-a.s.,
\end{equation*}
where $K^\mathbb{P}$ is a non-decreasing process null at $0$. As shown in the proof of Theorem 4.5 in Soner et al. \cite{STZ1}
one can find a universal $Z$ such that for each $\mathbb{P}\in\mathcal{P}_H$,
$$Z_t= Z^\mathbb{P}_t,\ 0\leq t\leq 1,\ \mathbb{P}-a.s..$$
Defining $Y=V$, from (\ref{estv}), $Y\in\mathbb{D}^\infty_H$. Similar to Lemma 3.1 in Possamai and Zhou \cite{PZ}, we deduce that $Z\in\mathbb{H}^2_H$.
Then, it suffices to verify that the family of non-decreasing processes $\{K^\mathbb{P}\}_{\mathbb{P}\in\mathcal{P}_H}$ satisfies the minimum condition (\ref{min}). For a fixed $\mathbb{P}\in\mathcal{P}_H$, $t\in[0, 1]$ and each $\mathbb{P}'\in\mathcal{P}_H(t^+, \mathbb{P})$, using the notations from the proof of Theorem \ref{rep}, we have
$$V_t-y^{\mathbb{P}'}_t(1,\xi)=\mathbb{E}^{\mathbb{Q}'}_t\bigg[\int^1_tM_sdK^{\mathbb{P}'}_s\bigg]\geq e^{-\mu}
\mathbb{E}^{\mathbb{Q}'}_t
[K^{\mathbb{P}'}_1-K^{\mathbb{P}'}_t],
$$
where 
$$\frac{d\mathbb{Q}'}{d\mathbb{P}'}\bigg|_{\mathcal{F}_t}=H^{\mathbb{P}'}_t:=\int^t_0\kappa^{\mathbb{P}'}_sdW^{\mathbb{P}'}_s,\ 0\leq t\leq 1,\ \mathbb{P}'-a.s.,$$
By Definition \ref{bmodef}, $H:=\{H^\mathbb{P'}\}_{\mathbb{P}'\in\mathcal{P}_H(t^+, \mathbb{P})}$ is a $BMO(\mathcal{P}_H(t^+, \mathbb{P}))$-martingale. Applying Lemma \ref{bmolemma2} to the family $\mathcal{E}(H)$, there exists a $p>1$ such that
$$\sup_{\mathbb{P}'\in\mathcal{P}_H(t^+, \mathbb{P})}\sup_{0\leq t\leq 1}\bigg|\bigg|\mathbb{E}^{\mathbb{P}'}_t\bigg[\bigg(
\frac{\mathcal{E}(H)_t}{\mathcal{E}(H)_1}\bigg)^\frac{1}{p-1}\bigg]\bigg|\bigg|_{L^\infty(\mathbb{P}')}\leq C_E,$$
then
\begin{align*}
\mathbb{E}^{\mathbb{P}'}_t[K^{\mathbb{P}'}_1-K^{\mathbb{P}'}_t]&\leq\mathbb{E}^{\mathbb{P}'}_t\bigg[\frac{\mathcal{E}(H^{\mathbb{P}'})_1}{\mathcal{E}(H^{\mathbb{P}'})_t}
(K^{\mathbb{P}'}_1-K^{\mathbb{P}'}_t)\bigg]^\frac{1}{2p-1}
\mathbb{E}^{\mathbb{P}'}_t\bigg[\bigg(\frac{\mathcal{E}(H^{\mathbb{P}'})_1}{\mathcal{E}(H^{\mathbb{P}'})_t}\bigg)^{-\frac{1}{2p-2}}
(K^{\mathbb{P}'}_1-K^{\mathbb{P}'}_t)\bigg]^\frac{2p-2}{2p-1}\\
&\leq\mathbb{E}^{\mathbb{Q}'}_t
[K^{\mathbb{P}'}_1-K^{\mathbb{P}'}_t]^\frac{1}{2p-1}
\mathbb{E}^{\mathbb{P}'}_t\bigg[\bigg(\frac{\mathcal{E}(H^{\mathbb{P}'})_1}{\mathcal{E}(H^{\mathbb{P}'})_t}\bigg)^{-\frac{1}{p-1}}\bigg]^\frac{p-1}{2p-1}
\mathbb{E}^{\mathbb{P}'}_t[(K^{\mathbb{P}'}_1-K^{\mathbb{P}'}_t)^2]^\frac{p-1}{2p-1}\\
&\leq C_E^\frac{p-1}{2p-1}C_2^\frac{p-1}{2p-1}e^{\frac{\mu}{2p-1}}(V_t-y^{\mathbb{P}'}_t(1,\xi))^\frac{1}{2p-1}.
\end{align*}
From (\ref{vp+}) and (\ref{vp++}), we obtain
\begin{align*}0\leq\sideset{}{^\mathbb{P}}\essinf_{\mathbb{P}\in\mathcal{P}_H(t^+, \mathbb{P})}\mathbb{E}^{\mathbb{P}'}_t[K^{\mathbb{P}'}_1-K^{\mathbb{P}'}_t]&\leq C(V_t-\sideset{}{^\mathbb{P}}\esssup_{\mathbb{P}\in\mathcal{P}_H(t^+, \mathbb{P})}y^{\mathbb{P}'}_t(1,\xi))^\frac{1}{2p-1}\\ &=C(V^+_t-\sideset{}{^\mathbb{P}}\esssup_{\mathbb{P}\in\mathcal{P}_H(t^+, \mathbb{P})}y^{\mathbb{P}'}_t(1,\xi))^\frac{1}{2p-1}=0,\ \mathbb{P}-a.s.,
\end{align*}
which is the desired result.
\hfill{}$\square$\\[6pt]
For each $\xi\in\mathcal{L}^\infty_H$, one can find a sequence $\{\xi^n\}_{n\in\mathbb{N}}\subset UC_b(\Omega)$, such that $||\xi^n-\xi||_{L^\infty_H}\rightarrow 0$. 
Thanks to a prior estimates, we have the following main result of the section.
\begin{theorem}\label{mresu}
Under (A1)-(A5) and for a given $\xi\in \mathcal{L}^\infty_H$, the 2BSDE (\ref{key}) has a unique solution $(Y, Z)\in \mathbb{D}^{\infty}_H\times\mathbb{H}^2_H$.
\end{theorem}
\begin{remark}
Recently, working under the Zermelo-Fraenkel set theory with axiom of choice (ZFC) and the continuum hypothesis (CH), Nutz \cite{N} has developed a way to define pathwisely a universal process that $\mathbb{P}$-a.s. coincides with the It\^o type stochastic integral of a predicable process $H$ with respect to a process $X$ that is a semimartingale under each $\mathbb{P}$ (for a continuous integrator $X$, $H$ needs only to be progressively measurable).
We notice that for each $\mathbb{P}\in\mathcal{P}_H$, the canonical process $B$ satisfies Assumption 2.1 in Nutz \cite{N}, then the stochastic integral $\int^T_tZ_sdB_s$ could be defined universally if we add ZFC and CH into our framework, but the aggregated process is only $\mathcal{F}^*$-adapted, where
$$
\mathcal{F}^*_t
:=\bigcap_{\mathbb{P}\in\mathcal{P}_H}
\mathcal{F}^+_t\vee\mathcal{N}^\mathbb{P}.
$$
On this occasion, $K$ could be a universal process in (\ref{key}). 
\end{remark}
\section{Application to finance}
\noindent In this section, we re-solve some robust utility maximization problems introduced by Matoussi et al. \cite{MPZ}.
\subsection{Statement of the problem} The problem under consideration in Matoussi et al. \cite{MPZ} is to maximize in a robust way the expected utility of the terminal value of a portfolio on a financial market with some uncertainty on the objective probability and to choose an optimal trading strategy to attain this optimal goal under some restrictions.\\[6pt]
This problem can be formulated into
\begin{align}\label{probe}
V(x):=\sup_{\pi\in\tilde{\mathcal{A}}}\inf_{\mathbb{P}\in\mathcal{P}}\mathbb{E}^{\mathbb{P}}[U(X^\pi_T-\xi)],
\end{align}
where $X^\pi_T$ is the terminal value of the wealth process associated with a strategy $\pi$ from a given set $\tilde{\mathcal{A}}$ of all admissible trading strategies, $\xi$ is a liability that matures at time $T$, $U$ denotes the utility function 
and $\mathcal{P}$ is a set of all possible probability measures. Without loss of generality, we always assume that $T=1$ in the sequel.\\[6pt]
In the present paper, we study the problem consists of non-dominated models i.e., the probability measures from the collection $\mathcal{P}$ could not be dominated by a finite measure.
Consistent with the setting for 2BSDE theory, we assume that $\mathcal{P}$ is a subset of the class $\bar{\mathcal{P}}_S$ (cf. Definition \ref{qs}), in which all the probability measures 
are mutually singular. 
\begin{definition}\label{p1}
In (\ref{probe}), let $\mathcal{P}=\tilde{\mathcal{P}}_H$ denote the collection of all those $\mathbb{P}\in \overline{\mathcal{P}}_S$ such that $$\underline{a}\leq\hat{a}_t\leq\overline{a}\ {\rm and}\ \hat{a}_t\in D_{F_t},\ \lambda\times \mathbb{P}-a.e.,
$$
for some $\underline{a}$, $\overline{a}\in\mathbb{S}_d^{>0}$ and each $(y,z)\in\mathbb{R}\times\mathbb{R}^d$.
\end{definition}
\noindent Adapted to this setting of $\tilde{\mathcal{P}}_H$, we shall change a little our settings for quadratic 2BSDEs, that is, (A3) and (A5) will be replaced by the following (A3') and (A5'):\\[6pt]
\textbf{(A3')} \textit{$F$ is continuous in $(y, z)$ and has a quadratic growth, i.e., for each $(\omega, t, y, z, a)\in\Omega\times[0, 1]\times\mathbb{R}\times\mathbb{R}^d\times D_{F_t}$,}
$$|F_t(\omega, y, z, a)|\leq \alpha(a)+\beta(a)|y|+\frac{\gamma}{2}|a^{1/2}z|^2,$$
where $\gamma$ is a strictly positive constant and $\alpha$, $\beta$ are non-negative deterministic functions satisfy that 
for some strictly positive constants $\overline{\alpha}$ and $\overline{\beta}$,
\begin{equation*}
\alpha(\hat{a}_t)\leq\overline{\alpha},\ {\rm and}\ \beta(\hat{a}_t)\leq\overline{\beta},\ \ 0\leq t\leq 1,\ \tilde{\mathcal{P}}_H-q.s..
\end{equation*}
\noindent\textbf{(A5')} $F$ is local Lipschitz in $z$, i.e., for each $(\omega, t, y, z, z', a)\in\Omega\times[0, 1]\times\mathbb{R}\times\mathbb{R}^d\times\mathbb{R}^d\times D_{F_t}$,
$$|F_t(\omega, y, z, a)-F_t(\omega, y, z', a)|\leq C (|a^{1/2}\phi(a)|+|a^{1/2}z|+|a^{1/2}z'|)|a^{1/2}(z-z')|
,$$
where $C$ is a strictly positive constant and $\bar{\phi}(a)=a^{1/2}\phi(a)$ satisfies that for some strictly positive
constant $\overline{\gamma}$,
$$|\bar{\phi}(\hat{a}_t)|=|\hat{a}^{1/2}_t\phi(\hat{a}_t)|\leq\overline{\gamma},\ 0\leq t\leq 1,\ \tilde{\mathcal{P}}_H-q.s..
$$
Repeating all the proof for the wellposedness of quadratic 2BSDEs in the last section, we can have the following theorem:
\begin{theorem}\label{mtho}
Under (A1)- (A2), (A3'), (A4) and (A5') and for a given $\xi\in \mathcal{L}^\infty_{\tilde{H}}$, the 2BSDE (\ref{key}) has a unique solution $(Y, Z)\in \mathbb{D}^{\infty}_H\times\mathbb{H}^2_H$.
\end{theorem}
\noindent We will have a detailed discussion for this kind of settings later in Subsection 5.4.\\[6pt]
\noindent The financial market consists of one bond with zero interest rate and $d$ stocks. The price process of the stocks is give by the following stochastic differential equations:
$$
dS^i_t=S^i_t(b^i_tdt+dB^i_t),\ 0\leq t\leq 1,\ i=1,\ldots,d,\ \tilde{\mathcal{P}}_H-q.s.,
$$
where $B$ is a $d$-dimensinal canonical process, $b^i$ is an $\mathbb{R}$-valued process that is uniformly bounded by a constant $M>0$ and is uniformly continuous in $\omega$ under the uniform norm $||\cdot||^\infty_1$, $i=1,2,\ldots,d$. By the definition of $\tilde{\mathcal{P}}_H$, for each $\mathbb{P}\in\tilde{\mathcal{P}}_H$, 
$$B_t=\int^t_0\hat{a}^{1/2}_sdW^\mathbb{P}_s,\ 0\leq t\leq 1,\ \mathbb{P}-a.s.,$$ 
where $\hat{a}^{1/2}$ plays in fact the role of volatility in (1) of Hu et al. \cite{HIM}. Thus, the difference of $\hat{a}^{1/2}$ under each $\mathbb{P}\in\tilde{\mathcal{P}}_H$ allows us to model the volatility uncertainty.\\[6pt]
\noindent In the following subsections, we study the problem (\ref{probe}) for two kinds of utility functions, the exponential and the power ones.
\subsection{Robust exponential utility maximization}
\noindent In this subsection, we consider the robust utility maximization problem (\ref{probe}) with an exponential utility function:
$$
U(x):=-\exp(-cx),\ c>0,\ x\in\mathbb{R}.
$$
In this case,
we denote $\pi=\{\pi_t\}_{0\leq t\leq 1}$ the trading strategy, which is a $d$-dimensional  $\mathcal{F}$-progressive measurable process. The $i$th component 
$\pi^i_t$ describes the amount of money invested in stock $i$ at time $t$, $i=1,\ldots,d$, then, for a given trading strategy $\pi$, the wealth process $X^\pi$ can be written as
$$
X^\pi_t=x+\sum^d_{i=1}\int^t_0\frac{\pi^i_s}{S^i_s}dS^i_s=x+\int^t_0\pi_s(dB_s+b_sds),\ 0\leq t\leq 1,\ \tilde{\mathcal{P}}_H-q.s..
$$
We now give the definition of the admissible trading strategies.
\begin{definition}\label{admexp}
Let $\tilde{C}$ be a closed set in $\mathbb{R}^d$. The set of admissible trading strategies $\tilde{\mathcal{A}}$ consists of all $d$-dimensional progressively measurable processes $\pi=\{\pi_t\}_{0\leq t\leq 1}$ that take values in $\tilde{C}$, $\lambda\otimes\mathcal{P}_H$-q.s., such that for each $\mathbb{P}\in\tilde{\mathcal{P}}_H$, $\int^1_0|\hat{a}^{1/2}_t\pi_t|^2dt<+\infty$, $\mathbb{P}$-a.s. and
$\{\exp(-cX^\pi_\tau)\}_{\tau\in\mathcal{T}^1_0}
$
is a $\mathbb{P}$-uniformly integrable family.
\end{definition}
\noindent Then, the utility maximization problem is equivalent to
\begin{equation}\label{uexp}
V(x):=\sup_{\pi\in\tilde{\mathcal{A}}}\inf_{\mathbb{P}\in\tilde{\mathcal{P}}_H}\mathbb{E}^{\mathbb{P}}\bigg[-\exp\bigg(-c\bigg(x+\int^1_0\pi_t(dB_t+b_tdt)-\xi\bigg)\bigg)\bigg].
\end{equation}
We can also consider a reduced utility maximization problem under each $\mathbb{P}\in\tilde{\mathcal{P}}_H$, which is introduced by Theorem 7 in Hu et al. \cite{HIM} and Theorem 4.1 in Morlais \cite{M}. Following these well known results, one can find a ${\pi^{\mathbb{P}}}^*\in\tilde{\mathcal{A}}^\mathbb{P}$ that solves the reduced utility maximization problem: 
\begin{equation}\label{suexp}
V^\mathbb{P}(x):=\sup_{\pi\in\tilde{\mathcal{A}^\mathbb{P}}}\mathbb{E}^{\mathbb{P}}\bigg[-\exp\bigg(-c\bigg(x+\int^1_0\pi_t(dB_t+b_tdt)-\xi\bigg)\bigg)\bigg],
\end{equation}
where $\tilde{\mathcal{A}}^\mathbb{P}$ is the collection of all admissible trading strategies given by Definition 1 in Hu et al. \cite{HIM} under $\mathbb{P}$ and thus, $\tilde{\mathcal{A}}\subset\tilde{\mathcal{A}}^\mathbb{P}$. 
It is evident that
$$
V(x)\leq\inf_{\mathbb{P}\in\tilde{\mathcal{P}}_H}V^\mathbb{P}(x).
$$
Therefore, the robust utility maximization problem (\ref{uexp}) is solved if one can find an optimal strategy $\pi^*$ such that 
$$V(x)=\inf_{\mathbb{P}\in\tilde{\mathcal{P}}_H}\mathbb{E}^{\mathbb{P}}\bigg[-\exp\bigg(-c\bigg(x+\int^1_0\pi^*_t(dB_t+b_tdt)-\xi\bigg)\bigg)\bigg]=\inf_{\mathbb{P}\in\tilde{\mathcal{P}}_H}V^\mathbb{P}(x).$$
In what follows, we give the theorem similar to Theorem 4.1 in Matoussi et al. \cite{MPZ} but without some additional condition on $\xi$ or on the border of $\tilde{C}$.
\begin{theorem}\label{exptheo}
Assume that $\xi\in\mathcal{L}^\infty_{\tilde{H}}$. The value function of the utility maximization problem (\ref{uexp}) is given by
$$
V(x)=-\exp(-c(x-Y_0)),
$$
where $Y_0$ is defined by the unique solution $(Y, Z)\in \tilde{\mathbb{D}}^\infty_{H}\times\tilde{\mathbb{H}}^2_{H}$ of the following 2BSDE:
\begin{equation}\label{exp2BSDE}
Y_t=\xi+\int^1_t\hat{F}_s(Z_s)ds-\int^1_tZ_sdB_s+K_1-K_t,\ 0\leq t\leq 1,\ \tilde{\mathcal{P}}_H-q.s.,
\end{equation}
where for each $(\omega, t, z, a)\in\Omega\times[0,1]\times\mathbb{R}^d\times\mathbb{S}^{>0}_d$,
\begin{equation}\label{expgen}
F_t(\omega, z, a):=\frac{c}{2}
\sideset{}{^2}\dist\bigg(a^{1/2}z
+\frac{1}{c}a^{-1/2}b_t(\omega),
a^{1/2}\tilde{C}\bigg)
-z^{\textnormal{\textbf{Tr}}}b_t(\omega)-\frac{1}{2c}|a^{-1/2}b_t(\omega)|^2.
\end{equation}
Moreover, there exists an optimal trading strategy $\pi^*\in\tilde{\mathcal{A}}$ with
\begin{equation}\label{optstra}
\hat{a}^{1/2}_t\pi^*_t\in\Pi_{\hat{a}^{1/2}_t\tilde{C}}\bigg(\hat{a}^{1/2}_tZ_t+\frac{1}{c}\hat{a}^{-1/2}_tb_t\bigg),\ \lambda\otimes\tilde{\mathcal{P}}_H-q.s.,
\end{equation}
where
$\Pi_{A}(r)$ denotes the collection of the elements in the closed set $A$ that realize the minimal distance to the point $r$.
\end{theorem}
\begin{remark} Some of the assumptions adopted by Theorem 4.1 in Matoussi et al. \cite{MPZ} are removed: our assumptions for the wellposedness of quadratic 2BSDEs does not concern the size of $\xi$, so we do not need to assume in addition that the liability $\xi$ is small enough in norm; on the other hand, we do not have any requirement on the regularity of the derivatives of the generator $\hat{F}$ and thus, 
the border of $\tilde{C}$ is no longer assumed to be a $\mathcal{C}^2$ curve.
It is evident that these two additional assumptions have limitations in real financial market: the one on $\xi$ is not practical;
the other one on the border of $\tilde{C}$ is often difficult to verify.
\end{remark}
\noindent\textbf{Sketch of the proof:} We prove this theorem by following procedures adopted by Matoussi et al. \cite{MPZ} but with some modifications, and we only give the sketch. \\[6pt]
\textbf{Step 1:} In this step, we show that the 2BSDE (\ref{exp2BSDE}) has a unique solution by verifying that the generator $F$ satisfies (A1)-(A2), (A3') and (A5'). Then, Theorem \ref{mtho} states that the 2BSDE (\ref{exp2BSDE}) admits a unique solution $(Y,Z)\in\tilde{\mathbb{D}}^\infty_H\times\tilde{\mathbb{H}}^2_H$.\\[3pt]
\begin{itemize}[leftmargin=*]
\item From that $b$ is uniform bounded and that $\tilde{C}$ is closed, we have, for each $(\omega, z)\in\Omega\times\mathbb{R}^d$, $D_{F_t(\omega, z)}=\mathbb{S}^{>0}_d$, which implies that (A1) is satisfied.\\[3pt]
\item Since $b$ is $\mathcal{F}$-progressive measurable and uniformly continuous in $\omega$ under the uniform norm, for each $(z, a)\in\mathbb{R}^d\times\mathbb{S}^{>0}_d$, $F(z, a)$ is $\mathcal{F}$-progressive measurable and uniformly continuous in $\omega$.\\[3pt]
\item For each $a\in\mathbb{S}^+_d$ that satisfies $\underline{a}\leq a \leq\overline{a}$, there exist a $\overline{K}>0$ that depends only on $\overline{a}$ and $\tilde{C}$ such that
\begin{equation*}
\inf\{|r|: r\in a^{1/2}\tilde{C}\}\leq \overline{K}, 
\end{equation*}
and another $\underline{K}>0$ that depends only on $\underline{a}$ and $M$, such that for each $\omega\in\Omega$, 
\begin{equation*}
|a^{-1/2}b_t(\omega)|^2\leq {\rm tr}(a^{-1})M^2=\underline{K}^2.
\end{equation*}
Then, for each $(t, z)\in[0,1]\times\mathbb{R}^d$,
\begin{equation}\label{est1}
\sideset{}{^2}\dist\bigg(a^{1/2}z
+\frac{1}{c}a^{-1/2}b_t,
a^{1/2}\tilde{C}\bigg)\leq 2|a^{1/2}z|^2+2\bigg(\frac{1}{c}|a^{-1/2}b_t|+\overline{K}\bigg)^2,
\end{equation}
from which we deduce 
$$|F_t(\omega, z, a)|\leq \bigg(2c\overline{K}^2+\frac{5+c}{2c}\underline{K}^2\bigg)+\bigg(\frac{1}{2}+c\bigg)|a^{1/2}z|^2.$$ 
That is to say (A3') is satisfied.\\[3pt]
\item
For each $(t, z^1, z^2)\in[0,1]\times\mathbb{R}^d\times\mathbb{R}^d$ and $a\in\mathbb{S}^{>0}_d$ that satisfies $\underline{a}\leq a \leq\overline{a}$,
\begin{align*}
F_t(\omega, z^1, a)-F_t(\omega, z^2, a)
=&\ \frac{c}{2}\bigg(\sideset{}{^2}\dist\bigg(a^{1/2}z^1
+\frac{1}{c}a^{-1/2}b_t,
a^{1/2}\tilde{C}\bigg)\\
-&\ \sideset{}{^2}\dist\bigg(a^{1/2}z^2
+\frac{1}{c}a^{-1/2}b_t,
a^{1/2}\tilde{C}\bigg)\bigg)-(z^1-z^2)^{\textnormal{\textbf{Tr}}}b_t.
\end{align*}
By the Lipschitz property of the distance function with respect to a closed set, we obtain the following inequality:
\begin{align*}
|F_t(\omega, z^1, a)-F_t(\omega, z^2, a)|&\leq \frac{c}{2}\bigg(\bigg(2\overline{K}+\frac{4}{c}\underline{K}\bigg)+|a^{1/2}z^1|+|a^{1/2}z^2|\bigg)|a^{1/2}(z^1-z^2)|,
\end{align*}
from which (A5') is satisfied.\\[3pt]
\end{itemize}

\noindent\textbf{Step 2:} We define, for each $\pi\in\tilde{\mathcal{A}}$,
\begin{equation}\label{rdef}
R^\pi_t=-\exp(-c(X^\pi_t-Y_t)),\ 0\leq t\leq 1,
\end{equation}
where $Y$ is the solution to 2BSDE (\ref{exp2BSDE}).
Then, we decompose $R^\pi$ into a product of two processes, i.e., $R^\pi=M^\pi A^\pi$, where for each $\mathbb{P}\in\tilde{\mathcal{P}}_H$,
\begin{align*}
M^\pi_t:=e^{-c(x-Y_0)}\exp\bigg(&-\int^t_0c(\pi_s-Z_s)dB_s\\&-\frac{1}{2}\int^t_0c^2|\hat{a}^{1/2}_s(\pi_s-Z_s)|^2ds-cK^\mathbb{P}_t\bigg),
\ 0\leq t\leq 1,\ \mathbb{P}-a.s.,
\end{align*}
and
$$
A^\pi_t:=-\exp\bigg(-\int^t_0\bigg(c\pi^{\textnormal{\textbf{Tr}}}_s b_s+c\hat{F}_s(Z_s)-\frac{1}{2}c^2|\hat{a}^{1/2}_s(\pi_s-Z_s)|^2\bigg)ds\bigg),\ 0\leq t\leq 1,\ \mathbb{P}-a.s..
$$
We rewrite $A^\pi$ into the following form,
\begin{align*}
A^\pi_t=-\exp\bigg(-\int^t_0\bigg(\frac{c^2}{2}\bigg|\hat{a}^{1/2}_s\pi_s-\bigg(\hat{a}^{1/2}_sZ_s&+\frac{1}{c}\hat{a}^{-1/2}_sb_s\bigg)\bigg|^2\\
&-cZ^{\textnormal{\textbf{Tr}}}_sb_s-\frac{1}{2}|\hat{a}^{-1/2}_sb_s|^2-c\hat{F}_s(Z_s)\bigg)ds\bigg).
\end{align*}
It is readily to observe that if $\pi=\pi^*$ that satisfies (\ref{optstra}), then 
$$A^{\pi^*}_t\equiv -1,\ 0\leq t\leq 1,\ \mathbb{P}-a.s..$$
%
Moreover, Lemma 11 in Hu et al. \cite{HIM} says that one can define such a $\pi^*$ that is 
$\mathcal{F}$-progressively measurable if $Z$ is $\mathcal{F}$-progressively measurable.\\[6pt]
In the previous section, we have already proved that $Z\in\tilde{\mathbb{H}}^2_{BMO(\tilde{\mathcal{P}}_H)}$. To show that $\pi^*-Z\in\tilde{\mathbb{H}}^2_{BMO(\tilde{\mathcal{P}}_H)}$, it suffices to verify that $\pi^*$ is also in $\tilde{\mathbb{H}}^2_{BMO(\tilde{\mathcal{P}}_H)}$. 
Applying triangle inequality to $|\hat{a}^{1/2}_t\pi^*_t|$ and recalling (\ref{est1}), we have, for each $t\in[0, 1]$,
\begin{align}\label{p3}
|\hat{a}^{1/2}_t\pi^*_t|&\leq\bigg|\hat{a}^{1/2}_tZ_t+\frac{1}{c}\hat{a}^{-1/2}_tb_t\bigg|+\bigg|\hat{a}^{1/2}_t\pi^*_t-\bigg(\hat{a}^{1/2}_tZ_t+\frac{1}{c}\hat{a}^{-1/2}_tb_t\bigg)\bigg|\notag\\
&\leq|\hat{a}^{1/2}_tZ_t|+\frac{1}{c}|\hat{a}^{-1/2}_tb_t|+
\dist\bigg(a^{1/2}_tZ_t+\frac{1}{c}\hat{a}^{-1/2}_tb_t, \hat{a}^{1/2}\tilde{C}\bigg)\\
&\leq 2|\hat{a}^{1/2}_tZ_t|+\frac{2}{c}\underline{K}+2\overline{K},\notag\ \mathcal{P}_H-q.s.,
\end{align}
which implies that $\pi^*$ is an element in $\tilde{\mathbb{H}}^2_{BMO(\tilde{\mathcal{P}}_H)}$.\\[6pt]
\noindent As $\pi^*\in\tilde{\mathbb{H}}^2_{BMO(\tilde{\mathcal{P}}_H)}$, for each $\mathbb{P}\in\tilde{\mathcal{P}}_H$, $\hat{a}^{1/2}\pi^*$ is a $BMO(\mathbb{P})$-martingale generator. By Remark 8 in Hu et al. \cite{HIM}, 
$\{\exp{-cX^\pi_\tau}\}_{\tau\in\mathcal{T}^1_0}$ is a $\mathbb{P}$-uniformly integrable family and it is easy to verify that $\mathbb{E}^\mathbb{P}[\int^1_0|\hat{a}^{1/2}_t\pi^*_t|^2dt]<+\infty$. Thus, $\pi^*\in \tilde{\mathcal{A}}$.\\[6pt]
\noindent\textbf{Step 3:} We now prove that for each $\mathbb{P}\in\tilde{\mathcal{P}}_H$,
\begin{equation}\label{est5}
\sideset{}{^\mathbb{P}}\esssup_{\mathbb{P}'\in\tilde{\mathcal{P}}_H(t, \mathbb{P})}\mathbb{E}^{\mathbb{P}'}_t[M^{\pi^*}_1]=M^{\pi^*}_t,\ 0\leq t\leq 1,\ \mathbb{P}-a.s.,
\end{equation}
so that
\begin{equation}\label{eq109}
\sideset{}{^\mathbb{P}}\essinf_{\mathbb{P}'\in\tilde{\mathcal{P}}_H(t, \mathbb{P})}\mathbb{E}^{\mathbb{P}'}_t[R^{\pi^*}_1]=R^{\pi^*}_t,\ 0\leq t\leq 1,\ \mathbb{P}-a.s..
\end{equation}
Since $-c(\pi^*-Z)$ is a $BMO(\tilde{\mathcal{P}}_H)$-martingale generator, under each $\mathbb{P}'\in\tilde{\mathcal{P}}_H(t, \mathbb{P})$, $\mathcal{E}(-c\int^\cdot_0(\pi^*_t-Z_t)dB_t)$ is an exponential martingale, and $M^{\pi^*}$ can be regard as a product of a martingale and a positive non-increasing process. Thus, it is easy to show that for each $\mathbb{P}\in\tilde{\mathcal{P}}_H$,
\begin{equation}\label{est4}
\sideset{}{^\mathbb{P}}\esssup_{\mathbb{P}'\in\tilde{\mathcal{P}}_H(t, \mathbb{P})}\mathbb{E}^{\mathbb{P}'}_t[M^{\pi^*}_1]\leq M^{\pi^*}_t,\ 0\leq t\leq 1,\ \mathbb{P}-a.s..
\end{equation}
To get the desired result, it suffices to prove the reverse inequality. Noticing that $M^{\pi^*}_1$ and $M^{\pi^*}_t$ are both positive, we can consider the ratio $\frac{M^{\pi^*}_1}{M^{\pi^*}_t}$. We calculate for each $t\in[0, 1]$ and $\mathbb{P}'\in\mathcal{P}_H(t, \mathbb{P})$,
\begin{align*}
\frac{M^{\pi^*}_1}{M^{\pi^*}_t}=\exp\bigg(&-\int^1_tc(\pi^*_s-Z_s)dB_s
\\&-\frac{1}{2}\int^1_tc^2|\hat{a}^{1/2}_s(\pi^*_s-Z_s)|^2ds-c(K^{\mathbb{P}'}_1-K^{\mathbb{P}'}_t)\bigg),\ \mathbb{P}'-a.s..
\end{align*}
Changing measure by
$$
\frac{d\mathbb{Q}'}{d\mathbb{P}'}\bigg|_{\mathcal{F}_t}=\mathcal{E}\bigg(-c\int^\cdot_0(\pi^*_s-Z_s)\hat{a}^{1/2}_sdW^{\mathbb{P}'}_s\bigg)_t,
$$
we have
$$\mathbb{E}^{\mathbb{P}'}_t\bigg[\frac{M^{\pi^*}_1}{M^{\pi^*}_t}\bigg]=\mathbb{E}^{\mathbb{Q}'}_t[\exp(-c(K^{\mathbb{P}'}_1-K^{\mathbb{P}'}_t))],\ \mathbb{P}'-a.s..
$$
By Jensen's inequality and the convexity of $\exp(-cx)$, we obtain
\begin{align*}
\sideset{}{^\mathbb{P}}\esssup_{\mathbb{P}'\in\tilde{\mathcal{P}}_H(t, \mathbb{P})}\mathbb{E}^{\mathbb{P}'}_t\bigg[\frac{M^{\pi^*}_1}{M^{\pi^*}_t}\bigg]
&=\sideset{}{^\mathbb{P}}\esssup_{\mathbb{P}'\in\tilde{\mathcal{P}}_H(t, \mathbb{P})}\mathbb{E}^{\mathbb{Q}'}_t[\exp(-c(K^{\mathbb{P}'}_1-K^{\mathbb{P}'}_t))]\\
&\geq\sideset{}{^\mathbb{P}}\esssup_{\mathbb{P}'\in\tilde{\mathcal{P}}_H(t, \mathbb{P})}\exp(-c\mathbb{E}^{\mathbb{Q}'}_t[K^{\mathbb{P}'}_1-K^{\mathbb{P}'}_t])\\
&\geq \exp(-c\sideset{}{^\mathbb{P}}\essinf_{\mathbb{P}'\in\tilde{\mathcal{P}}_H(t, \mathbb{P})}\mathbb{E}^{\mathbb{Q}'}_t[K^{\mathbb{P}'}_1-K^{\mathbb{P}'}_t]).
\end{align*}
Similar to (\ref{est2}), we know, for some $p$, $q>1$ that satisfy $1/{p}+{1}/{q}=1$,
$$
\sideset{}{^\mathbb{P}}\essinf_{\mathbb{P}'\in\tilde{\mathcal{P}}_H(t, \mathbb{P})}\mathbb{E}^{\mathbb{Q}'}_t[K^{\mathbb{P}'}_1-K^{\mathbb{P}'}_t]
\leq C_{RH}^{1/q}C^{1/2p}_{2p-1}\sideset{}{^\mathbb{P}}\essinf_{\mathbb{P}'\in\tilde{\mathcal{P}}_H(t, \mathbb{P})}\mathbb{E}^{\mathbb{P}'}_t[K^{\mathbb{P}'}_1-K^{\mathbb{P}'}_t]=0,
$$
where $C_{RH}$ is the constant in Lemma \ref{BMOP} and $C_{2p-1}$ is from (\ref{18}). The inequality above implies that
\begin{equation}\label{est3}
\sideset{}{^\mathbb{P}}\esssup_{\mathbb{P}'\in\tilde{\mathcal{P}}_H(t, \mathbb{P})}\mathbb{E}^{\mathbb{P}'}_t\bigg[\frac{M^{\pi^*}_1}{M^{\pi^*}_t}\bigg]\geq 1.
\end{equation}
Then, (\ref{est5}) comes after (\ref{est4}) and (\ref{est3}).\\[6pt]
\noindent\textbf{Step 4:} Under each $\mathbb{P}\in\tilde{\mathcal{P}}_H$, the canonical process $B$ is a $\mathbb{P}$-martingale and $\hat{F}_t(z)$ is in fact (2.6) in Morlais \cite{M}. Thus, the value function of the reduced utility maximization problem is given by 
$$
V^\mathbb{P}(x)=-\exp(-c(x-Y^\mathbb{P}_0)),
$$
where $Y^\mathbb{P}_0$ is defined by the unique solution $(Y^\mathbb{P}, Z^\mathbb{P})\in D^\infty({\mathbb{P}})\times H^2(\mathbb{P})$ of the following BSDE:
\begin{equation}\label{eq108}Y^\mathbb{P}_t=\xi+\int^1_t\hat{F}_s(Z_s)ds-\int^1_tZ_sdB_s,\ 0\leq t\leq 1,\ \mathbb{P}-a.s.. 
\end{equation}
By Theorem 3.2, we have
$$Y_0=\sup_{\mathbb{P}\in\tilde{\mathcal{P}}_H} Y^\mathbb{P}_0.$$
From (\ref{eq109}) and (\ref{eq108}), it holds true that
\begin{align*}
\inf_{\mathbb{P}\in\tilde{
\mathcal{P}}_H}
\mathbb{E}^\mathbb{P}[-\exp(-c(X^{\pi*}_t-\xi))]&=\inf_{\mathbb{P}
\in\tilde{\mathcal{P}}_H}
\mathbb{E}^\mathbb{P}[R^{
\pi^*}_1]=R^{\pi^*}_0\\&=-\exp(-c(x-Y_0))=\inf_{\mathbb{P}\in\tilde{
\mathcal{P}}_H} -\exp(-c(x-Y^\mathbb{P}_0)),
\end{align*}
which implies that $\pi^*$ is the optimal strategy. We complete the proof.\hfill$\square$
\begin{remark}
In fact, we adopt a weaker assumption on the admissible strategy than the one in Theorem 4.1 in Matoussi et al. \cite{MPZ}. We only assume that $\pi$ is an admissible strategy defined by Hu et al. \cite{HIM} and Morlais \cite{M} under each $\mathbb{P}\in\tilde{\mathcal{P}}_H$, i.e., $$\tilde{\mathcal{A}}=\bigcap_{\mathbb{P}\in \tilde{\mathcal{P}}_H}\tilde{\mathcal{A}}^\mathbb{P},$$
while Matoussi et al. \cite{MPZ} assumed that $\pi\in
\tilde{\mathbb{H}}^2_{BMO(
\tilde{\mathcal{P}}_H)}$. Under this stronger assumption, all $R^\pi$ satisfies the minimal condition (\ref{eq109}) and they verified that $\pi^*$ is optimal only $A^\pi\leq A^* \equiv -1$, for all $\pi$ is admissible. In our present paper, we justify that $\pi^*$ is optimal for this larger set of admissible strategies by a min-max property as we showed in Step 4, which is regardless of whether the admissible strategy other than the optimal one is an $BMO(\tilde{\mathcal{P}}_H)$-martingale generator. Although we have still proved that $\pi^*\in
\tilde{\mathbb{H}}^2_{BMO(
\tilde{\mathcal{P}}_H)}$, this result is more general.
\end{remark}
\subsection{Robust power utility maximization} In this subsection, we redo the problem (\ref{probe}) with a power utility function:
$$
U(x):=\frac{1}{\gamma}x^\gamma,\ \gamma<1,\ x\in\mathbb{R}.
$$
In this case, a $d$-dimensional $\mathcal{F}$-progressively measuable process $\{\rho_t\}_{0\leq t \leq 1}$ denotes the trading strategy, whose component $\rho^i_t$ describes the proportion of money invested in stock $i$ at time $t$, $0\leq t\leq 1$, $i=1, 2, \ldots, d$, then, for a given trading strategy $\rho$, the wealth process $X^\rho$ can be written as
\begin{equation}\label{wealth2}
X^\rho_t=x+\sum^d_{i=1}\int^t_0\frac{X^\rho_s\rho^i_s}{S^i_s}dS^i_s=x+\int^t_0 X^\rho_s\rho_s(dB_s+b_sds),\ 0\leq t\leq 1,\ \tilde{\mathcal{P}}_H-q.s.,
\end{equation}
where the initial capital $x$ is positive. One can find an $X^\rho$ defined by 
$$
X^\rho_t:=x\mathcal{E}\bigg(\int^\cdot_0\rho_s(dB_s+b_sds)\bigg)_t,\ 0\leq t\leq 1,
$$ 
which is the unique solution of (\ref{wealth2}) under each $\mathbb{P}\in\tilde{\mathcal{P}}_H$.
\begin{definition}\label{admpow}
Let $\tilde{C}$ be a closed set in $\mathbb{R}^d$. The set of admissible trading strategies $\tilde{\mathcal{A}}$ consists of all $d$-dimensional progressively measurable processes $\rho=\{\rho_t\}_{0\leq t\leq 1}$ that take values in $\tilde{C}$, $\lambda\otimes\tilde{\mathcal{P}}_H$-q.s. and for each $\mathbb{P}\in\tilde{\mathcal{P}}_H$, $\int^1_0|\hat{a}^{1/2}_t\rho_t|^2dt<+\infty$, $\mathbb{P}$-a.s..
\end{definition}
\noindent For each $\mathbb{P}\in\tilde{\mathcal{P}}_H$, we define a probability measure $\mathbb{Q}\ll\mathbb{P}$ by
$$
\frac{d\mathbb{Q}}{d\mathbb{P}}
\bigg|_{\mathcal{F}_t}
=\mathcal{E}\bigg(-\int^\cdot_0b^\textnormal{\textbf{Tr}}_s\hat{a}^{1/2}_sdW^\mathbb{P}_s\bigg)_t,
$$
then, by the definition above, for each $\rho\in\tilde{\mathcal{A}}$, $X^\rho$ is a $\mathbb{Q}$-local martingale bounded from below. Thus, $X^\rho$ is a $\mathbb{Q}$-supermartingale. Since $\mathbb{Q}\ll\mathbb{P}$, the strategy $\rho$ is free of arbitrage under $\mathbb{P}$.\\[6pt]
\noindent We suppose that the investor has no liability, i.e., $\xi=0$, then the maximization problem is equivalent to 
\begin{equation}\label{upow}
V(x):=\frac{1}{\gamma}x^\gamma\sup_{\rho\in\tilde{\mathcal{A}}}\inf_{\mathbb{P}\in\tilde{\mathcal{P}}_H}\mathbb{E}^\mathbb{P}\bigg[\exp\bigg(\gamma\int^1_0\rho_s(dB_s+b_sds)-\frac{\gamma}{2}\int^1_0|\hat{a}^{1/2}_s\rho_s|^2ds\bigg)\bigg].
\end{equation}
Similar to that in the last subsection, we have the following theorem:
\begin{theorem}\label{powtheo}
The value function of the utility maximization problem (\ref{upow}) is given by
$$
V(x)=\frac{1}{\gamma}x^\gamma\exp(Y_0),
$$
where $Y_0$ is defined by the unique solution $(Y, Z)\in \tilde{\mathbb{D}}^\infty_H\times\tilde{\mathbb{H}}^2_H$ of the following 2BSDE:
\begin{equation}\label{pow2BSDE}
Y_t=0+\int^1_t\hat{F}_s(Z_s)ds-\int^1_tZ_sdB_s+K_1-K_t,\ 0\leq t\leq 1,\ \tilde{\mathcal{P}}_H-q.s.,
\end{equation}
where for each
$(\omega, t, z, a)\in\Omega\times[0,1]\times\mathbb{R}^d\times\mathbb{S}^{>0}_d$,
\begin{align}\label{powgen}
F_t(\omega, z, a):=-\frac{\gamma(1-\gamma)}{2}
\sideset{}{^2}\dist\bigg(\frac{1}{1-\gamma}(a^{1/2}z
&+a^{-1/2}b_t(\omega)),
a^{1/2}_t\tilde{C}\bigg)\\
&+\frac{\gamma|a^{1/2}z
+a^{-1/2}b_t(\omega)|^2}{2(1-\gamma)}+\frac{1}{2}|a^{1/2}z|^2.\notag
\end{align}
Moreover, there exists an optimal trading strategy $\rho^*\in\tilde{\mathcal{A}}$ with
\begin{align}\label{optpow}
\hat{a}^{1/2}_t\rho^*_t(\omega)\in\Pi_{\hat{a}^{1/2}_t\tilde{C}}\bigg(\frac{1}{1-\gamma}(a^{1/2}_tz
+a^{-1/2}_tb_t(\omega))\bigg),\ 0\leq t\leq 1,\ \tilde{\mathcal{P}}_H-q.s.,
\end{align}
where
$\Pi_{A}(r)$ denotes the collection of the elements in the closed set $A$ that realize the minimal distance to the point $r$.
\end{theorem}
\noindent\textbf{Sketch of the proof:} Following similar procedures in the proof of Theorem \ref{exptheo}, we verify that the generator $F$ in 2BSDE (\ref{pow2BSDE}) satisfies (A1)-(A2), (A3') and (A5') and define a family of processes $\{R^\rho\}_{\rho\in\tilde{\mathcal{A}}}$ by
$$
R^\rho_t:=\frac{1}{\gamma}x^\gamma\exp\bigg(\gamma\int^t_0\rho_s(dB_s+b_sds)-\frac{\gamma}{2}\int^t_0|\hat{a}^{1/2}_s\rho_s|^2ds+Y_t\bigg),\ 0\leq t\leq 1,
$$
\noindent such that for each $\mathbb{P}\in\tilde{\mathcal{P}}_H$,
\begin{itemize}[leftmargin=*]
\item $R^\rho_0$ is a constant indepent of $\rho$;
\item $R^\rho_1=\frac{1}{\gamma}(X^\rho_1)^\gamma$, for each $\rho\in\tilde{\mathcal{A}}$.
\end{itemize}
Then, we rewrite $R^\rho$ under each $\mathbb{P}\in\tilde{\mathcal{P}}_H$ as follows:
$$
R^\rho_t=\frac{1}{\gamma}x^\gamma\exp(Y_0)\mathcal{E}\bigg(\int^\cdot_0(\gamma\rho_s+Z_s)dB_s\bigg)_te^{-K^\mathbb{P}_t}\exp\bigg(\int^t_0\nu_sds\bigg),\ 0\leq t\leq 1,\ \mathbb{P}-a.s.,
$$ 
where
\begin{align*}
\nu_t:=-\frac{\gamma(1-\gamma)}{2}
\bigg|\hat{a}^{1/2}_t\rho_t&-\frac{1}{1-\gamma}(\hat{a}^{1/2}_tZ_t
+\hat{a}^{-1/2}_tb_t)\bigg|^2\\&+\frac{\gamma|\hat{a}^{1/2}_tZ_t
+\hat{a}^{-1/2}_tb_t|^2}{2(1-\gamma)}+\frac{1}{2}|\hat{a}^{1/2}_tZ_t|^2-\hat{F}_t(Z_t).
\end{align*}
Similar to (\ref{eq109}), we could find an optimal strategy $\rho^*$ such that for each $\mathbb{P}\in\tilde{\mathcal{P}}_H$,
$$\nu^{\rho^*}_t\equiv 0,\ 0\leq t\leq 1,\ \mathbb{P}-a.s.$$
and thus,
\begin{equation}\label{lastass}
\sideset{}{^\mathbb{P}}\essinf_{\mathbb{P}'\in\tilde{\mathcal{P}}_H(t, \mathbb{P})}\mathbb{E}^{\mathbb{P}'}_t[R^{\rho^*}_1]= R^{\rho^*}_t,\ 0\leq t\leq 1,\ \mathbb{P}-a.s..
\end{equation}
The desired result comes after (\ref{lastass}) and the min-max property.\hfill{}$\square$
\begin{remark}
In Matoussi et al. \cite{MPZ}, only the case that $\gamma<0$ was considered. According to their assumption that $\tilde{C}$ contains $0$, we calculate
$$\hat{F}^0_t=-\frac{\gamma}{2(1-\gamma)}|a^{-1/2}_tb_t|^2,$$
where $-\frac{\gamma}{2(1-\gamma)}$ is dominated by $\frac{1}{2}$ when $\gamma<0$ and so that they can give a uniform assumption on $b$, which is regardless of $\gamma$, to make sure that $F^0$ is small enough.  
\end{remark}
\subsection{Some remarks on the class of probability measures and the assumptions}
We have already seen that the 2BSDEs (\ref{exp2BSDE}) and (\ref{pow2BSDE}) are discussed under some new settings, where $\mathcal{P}_H$ was changed into $\tilde{\mathcal{P}}_H$; (A3) and (A5) were changed into (A3') and (A5'). In what follows, we would like to discuss more about these conditions and class and probability measures. \\[6pt]
Since these weakened conditions shall be related to some given series of probability measure classes $\{\mathcal{P}^t_H\}_{t\in[0,1]}$ of probability measures, we first give the following definition:
\begin{definition}\label{consis} We say a series of probability measure classes $\{\mathcal{P}^t_H\}_{t\in[0,1]}$ is consistent if the following points are satisfied (we note $\mathcal{P}^0_H=\mathcal{P}_H$.):
\begin{itemize}[leftmargin=*]
\item For each $\mathbb{P}\in\mathcal{P}_H$, for $\mathbb{P}$-a.e. $\omega\in\Omega$ and each $\tau\in\mathcal{T}^1_0$, $\mathbb{P}^{\tau, \omega}\in\mathcal{P}^{\tau(\omega)}_H$;
\item For each $\tau\in\mathcal{T}^1_0$, $A\in\mathcal{F}_\tau$, $\mathbb{P}\in\mathcal{P}_H$ and $\hat{\mathbb{P}}^\tau\in\mathcal{P}^\tau_H$, $\mathbb{P}\otimes^A_\tau\hat{\mathbb{P}}^\tau\in\mathcal{P}_H$, where for each $E\subset\Omega$, 
$$
\mathbb{P}\otimes^A_\tau\hat{\mathbb{P}}^\tau(E):= \mathbb{E}^\mathbb{P}
[\mathbb{E}^{\hat{\mathbb{P}}^\tau}[({\bf 1}_E)^{\tau, \omega}]{\bf 1}_{A}]+\mathbb{P}(E\cap A^c). 
$$
\end{itemize}
\end{definition}
\noindent In the 2BSDE framework, the series of classes defined by Definition \ref{qs} is consistent, since the first point is guaranteed by  Lemma 4.1 in Soner et al.  \cite{STZ1} and the second one is in fact the reduced version ($n=1$) of the statement (4.19) in Soner et al. \cite{STZ1}. These two properties play an important role in our proof of the dynamic programming principle (cf. Proposition \ref{repv}).\\[6pt]
In what follows, we verify that the series of classes defined by Definition \ref{p1} is consistent.
In this case, $\tilde{\mathcal{P}}^t_H$ consists of 
all those $\mathbb{P}\in \overline{\mathcal{P}}^t_S$ such that $$\underline{a}\leq\hat{a}^t_s\leq\overline{a}\ {\rm and}\ \hat{a}^t_s\in D_{F_s},\ \lambda\times \mathbb{P}^t-a.e.,
$$
for some $\underline{a}$, $\overline{a}\in\mathbb{S}_d^{>0}$ and each $(y,z)\in\mathbb{R}\times\mathbb{R}^d$. Since $\tilde{\mathcal{P}}_H\subset \overline{\mathcal{P}}_S$, by Lemma 4.1 in Soner et al. \cite{STZ1},  for a given $\mathbb{P}\in\tilde{\mathcal{P}}_H$ and $\mathbb{P}$-a.e. $\omega\in\Omega$, $\mathbb{P}^{\tau, \omega}\in\overline{\mathcal{P}}^{\tau(\omega)}_S$ and
$$
\underline{a}\leq \hat{a}^{\tau(\omega)}_t(\tilde{\omega})=\hat{a}^{\tau, \omega}_t (\tilde{\omega})=\hat{a}_{t}(\omega\otimes^\tau\tilde{\omega})\leq \overline{a},\ \lambda\times\mathbb{P}^{\tau, \omega}-a.e..
$$
On the other hand, the proof of statement (4.19) in Soner et al. \cite{STZ1} showed that $\mathbb{P}\otimes^A_\tau\hat{\mathbb{P}}^\tau\in\overline{\mathcal{P}}_S$. Defining $\tilde{\mathbb{P}}:=\mathbb{P}
\otimes^A_\tau
\hat{\mathbb{P}}^\tau$, it suffices to verify that 
\begin{equation}\label{conl}
\underline{a}\leq \hat{a}_t\leq \overline{a},\ \lambda\times\tilde{\mathbb{P}}-a.e..
\end{equation}
We calculate
\begin{align*}
\int^1_0\mathbb{E}^{\tilde{\mathbb{P}}}[{\bf 1}_{\{\hat{a}_t\notin[\underline{a}, \overline{a}]\}}]dt
=\int^1_0(\mathbb{E}^\mathbb{P}
[\mathbb{E}^{\hat{\mathbb{P}}^\tau}[({\bf 1}_{\{\hat{a}_t\notin[\underline{a}, \overline{a}]\}})^{\tau, \omega}]{\bf 1}_{A}]
+
\mathbb{E}^\mathbb{P}[{\bf 1}_{\{\{\hat{a}_t\notin[\underline{a}, \overline{a}]\}\cap A^c\}}])dt,
\end{align*}
where 
$$
\mathbb{E}^{\hat{\mathbb{P}}^\tau}[({\bf 1}_{\{\hat{a}_t\notin[\underline{a}, \overline{a}]\}})^{\tau, \omega}]{\bf 1}_{A}(\omega)=\left\{
\begin{array}{l@{\quad, \quad}l}
\mathbb{E}^{\hat{\mathbb{P}}^\tau}[{\bf 1}_{\{\hat{a}^\tau_t\notin[\underline{a}, \overline{a}]\}}(\tilde{\omega})]=0 & \omega\in A,\ t\geq \tau(\omega);\\[3pt]
{\bf 1}_{\{\hat{a}_t\notin[\underline{a}, \overline{a}]\}} (\omega) & \omega\in A,\ t< \tau(\omega);\\[3pt]
0 & otherwise.
\end{array}\right.
$$
Thus, 
$$
\int^1_0\mathbb{E}^{\tilde{\mathbb{P}}}[{\bf 1}_{\{\hat{a}_t\notin[\underline{a}, \overline{a}]\}}]dt\leq \int^1_0\mathbb{E}^{{\mathbb{P}}}[{\bf 1}_{\{\hat{a}_t\notin[\underline{a}, \overline{a}]\}}]dt=0,
$$
which implies (\ref{conl}).
%
%
\begin{remark}\label{rem1} Suppose that a consistent series of probability measure classes $\{\mathcal{P}^t_H\}_{t\in[0, 1]}\subset \bar{\mathcal{P}}_S$ is given (not limited to the form defined by Definition \ref{qs} and \ref{p1}), then (A3) can be even weakened to the following form, which is similar to (H1) in Morlais \cite{M} for quadratic BSDEs:\\[6pt]
\textbf{(A3'')} \textit{$F$ is continuous in $(y, z)$ and has a quadratic growth in $z$, i.e., for each $(\omega, t, y, z, a)\in\Omega\times[0, 1]\times\mathbb{R}\times\mathbb{R}^d\times D_{F_t}$,}
\begin{equation}\label{efer}|F_t(\omega, y, z, a)|\leq \alpha_t(a)+\beta_t(a)|y|+\frac{\gamma}{2}|a^{1/2}z|^2,
\end{equation}
where $\gamma$ is a strictly positive constant and $\alpha$, $\beta$ satisfy that
\begin{itemize}[leftmargin=*]
\item For each $a\in\mathbb{S}^{>0}_d$, $\alpha(a)$, $\beta(a)$ are nonnegative
$\mathcal{F}$-progressive measurable processes;
\item For some $\overline{\alpha}$, $\overline{\beta}$ which are strictly positive constants,
\begin{equation*}\label{weakve}
\int^1_0\alpha_t(\hat{a}_t)dt\leq\overline{\alpha}\ {\rm and}\ \int^1_0\beta_t(\hat{a}_t)dt\leq\overline{\beta},\ \mathcal{P}_H-q.s.;
\end{equation*}
\item For each $(\omega, t)\in\Omega\times(0, 1]$ and $\mathbb{P}^t\in\mathcal{P}^t_H$, 
\begin{equation*}\label{weakver}
\int^1_t\alpha^{t, \omega}_s(\hat{a}^t_s)ds\leq\overline{\alpha}\ {\rm and}\ \int^1_t\beta^{t, \omega}_s(\hat{a}^t_s)ds\leq\overline{\beta},\ \mathbb{P}^t-a.s.,
\end{equation*}
\noindent where $\overline{\alpha}$, $\overline{\beta}$ are the same as above.\\
\end{itemize}
We recall (\ref{v}) that for each $(\omega, t)\in\Omega\times[0, 1]$, $V_t(\omega)$ concerns the solutions of the $(t, \omega)$-shifted quadratic BSDEs under all $\mathbb{P}^t\in\mathcal{P}^t_H$. 
Therefore, for each $t\in[0, 1]$, at least $\mathcal{P}_H$-q.s. $\omega\in\Omega$, $(t, \omega)$ shifted generator should satisfy (H1) in Morlais \cite{M} (or similar conditions for quadratic BSDEs) under each $\mathbb{P}^t\in\mathcal{P}^t_H$ to ensure the existence of these solutions. We notice that the orignal condition (A3) is posed pathwisely, that is, it holds for all $(\omega, t)\in\Omega\times[0, 1]$, whereas (A3'') also involves pathwise settings for each $(\omega, t)\in\Omega\times[0, 1]$. Therefore, (\ref{v}) can be well defined under these two conditions.\\[6pt]
A natural question arises: if (\ref{gere}) and (\ref{efer}) can be written in a $\mathcal{P}_H$-q.s. version; if the third point of (A3'') can be removed?\\[6pt]
We consider the first question: suppose that for all $(t, y, z, a)\in[0, 1]\times\mathbb{R}\times\mathbb{R}^d$,
\begin{equation*}
|\hat{F}_t(y, z)|\leq \alpha+\beta|y|+\frac{\gamma}{2}|\hat{a}^{1/2}_tz|^2,\ \mathcal{P}_H-q.s..
\end{equation*}
Fixing an $\mathbb{P}^t\in\mathcal{P}^t_H$, we can choose an arbitrage $\mathbb{P}\in\mathcal{P}_H$ and construct a concatenation probability $\hat{\mathbb{P}}:=\mathbb{P}\otimes^\Omega_t\mathbb{P}^t$. Since $\{\mathcal{P}^t_H\}_{t\in[0, 1]}$ is consistent, $\hat{\mathbb{P}}\in\mathcal{P}_H$, $\mathbb{P}|_{\mathcal{F}_t}=\hat{\mathbb{P}}|_{\mathcal{F}_t}$ and for each  $\omega\in\Omega$, $\hat{\mathbb{P}}^{t, \omega}=\mathbb{P}^t$. Thus, we have for $\mathbb{P}$-a.s, $\omega\in\Omega$ and all $(s, y, z, a)\in[t, 1]\times\mathbb{R}\times\mathbb{R}^d$,
\begin{align}\label{feret}
|\hat{F}^{t, \omega}_s(y, z)|&=|F_s(\omega\otimes_t\tilde{\omega}, y, z, \hat{a}_s(\omega\otimes_t\tilde{\omega}))
\\ &\leq\alpha+\beta|y|+\frac{\gamma}{2}|\hat{a}^{1/2}_s(\omega\otimes_t\tilde{\omega})z|^2=\alpha+\beta|y|+\frac{\gamma}{2}|(\hat{a}^t_s)^{1/2}(\tilde{\omega})z|^2,\ \mathbb{P}^t-a.s..\notag
\end{align}
Since $\mathbb{P}$ is arbitrage, we can deduce that for $\mathcal{P}_H$-q.s. $\omega\in\Omega$, (\ref{feret}) is satisfied.
In other words, defining for each $\mathbb{P}^t\in\mathcal{P}^t_H$ a set:
$$
E^{\mathbb{P}^t}:=\{\omega:\ 
\hat{F}^{t, \omega}_s(y, z)\ satisfies\ (\ref{feret}),\ \mathbb{P}^t-a.s.\},
$$
we have $\mathbb{P}(E^{\mathbb{P}^t})=1$ for all $\mathbb{P}\in\mathcal{P}_H$. At the end of the day, we still have no idea about $\mathbb{P}(\cap_{\mathbb{P}^t\in\mathcal{P}^t_H}E^{\mathbb{P}^t})$, since it is a probability of an intersection of non-countable sets. Therefore, the answer to the first question is negative.\\[6pt] 
For the same reason, the answer to the second question is negative either, unless we could find an $\alpha$ such that for each $a\in\mathbb{S}^{>0}_d$, $\alpha^{t, \omega}(a)$ is independent of $\omega$, i.e., $\alpha^{t, \omega}_s(a)\equiv \alpha^{t}_s(a)$.  In such case, if we only assume the second point and define
$$E^{\mathbb{P}^t}:=\bigg\{\omega:
\int^1_t\alpha^{t, \omega}_s(\hat{a}^t_s)ds\leq\overline{\alpha},\ \mathbb{P}^t-a.s.
\bigg\},
$$
then $E^{\mathbb{P}^t}=\Omega$ for all ${\mathbb{P}^t}\in\mathcal{P}^t_H$, which implies the third point in (A3'').
As we have shown in 
(A3'), a special case of such $\alpha$ is that
for each $a\in\mathbb{S}^{>0}_d$, $\alpha_\cdot(a)$ is a deterministic function in $t$.
\end{remark}
\begin{remark}
Corresponding to (v) of Assumption 2.2 in Possamai and Zhou \cite{PZ}, (A5) can be weakened to the following form:\\[6pt] 
\noindent\textbf{(A5'')} $F$ is local Lipschitz in $z$, i.e., for each $(\omega, t, y, z, z', a)\in\Omega\times[0, 1]\times\mathbb{R}\times\mathbb{R}^d\times\mathbb{R}^d\times D_{F_t}$,
$$|F_t(\omega, y, z, a)-F_t(\omega, y, z', a)|\leq C (|a^{1/2}\phi_t(a)|+|a^{1/2}z|+|a^{1/2}z'|)|a^{1/2}(z-z')|
,$$
where $C$ is a strictly positive constant and $\phi$ satisfies that 
\begin{itemize}[leftmargin=*]
\item For each $a\in\mathbb{S}^{>0}_d$, $\phi(a)$ is an $\mathcal{F}$-progressively measurable process;
\item $\phi(\hat{a})$ is a $BMO(\mathcal{P}_H)$-martingale generator;
\item For each $(\omega, t)\in\Omega\times(0, 1]$, $\phi^{t, \omega}(\hat{a}^t)$ is a $BMO(\mathbb{P}^t)$-martingale generator under each $\mathbb{P}^t\in\mathcal{P}^t_H$.\\
\end{itemize}
Based on the argument in remark \ref{rem1}, only having that $\phi$ is a $BMO(\mathcal{P}_H)$-martingale generator, we have no idea weather for some $(\omega, t)\in\Omega\times[0, 1]$, $\phi^{t, \omega}$ is a $BMO(\mathbb{P}^t)$-martingale generator under all $\mathbb{P}^t\in\mathcal{P}^t_H$, unless for each $a\in\mathbb{S}^{>0}_d$, $\phi(a)$ is independent of $\omega$. Thus, the third point in (A5'') is necessary.
We would like to point out that (v) in Assumption 2.2 in Possamai and Zhou \cite{PZ} is ambiguous, which may cause some slight problems for their setting of $V_t(\omega)$ and for the proof of Lemma 5.1 in that paper. 
\end{remark}
\noindent Taking the 2BSDE (\ref{exp2BSDE}) as an example, we explain these settings ((A3') and (A5') are special cases of (A3'') and (A5''), respectively). We observe that the generator (\ref{expgen}) satisfies the quadratic condition (A3'') for
$$
\alpha_t(a):=2c\inf\{|r|^2: r\in a^{1/2}\tilde{C}\}+\frac{5+c}{2c} {\rm tr}(a^{-1})M^2,\ a\in\mathbb{S}^{>0}_d, 
$$
in which $\alpha(a)$ is a deterministic function. In general, $\inf\{|r|^2: r\in a^{1/2}\tilde{C}\}$ and $|a^{-1/2}b|^2$ could be unbounded, so that (A3) is no longer satisfied. If we choose 
$\overline{\alpha}=2c\overline{K}^2+\frac{5+c}{2c}\underline{K}^2$, then (A3'') is satisfied.\\[6pt]
Similarly, the generator (\ref{expgen}) satisfies no longer and (A5). We define
$$\phi_t(a):=2\inf\{|r|: r\in a^{1/2}\tilde{C}\} +\frac{4}{c} ({\rm tr}(a^{-1}))^{1/2}M,\ a\in\mathbb{S}^{>0}_d,$$ which is bounded by $2\overline{K}+\frac{4}{c}\underline{K}$ when $a$ is replaced by $\hat{a}$ (or $\hat{a}^t$, respectively), $\tilde{\mathcal{P}}_H$ (or $\tilde{\mathcal{P}}^t_H$, respectively)-q.s.. By Definition \ref{p1}, we know that a constant process is a $BMO(\tilde{\mathcal{P}}_H)$ (or $BMO(\tilde{\mathcal{P}}^t_H)$, respectively)-martingale generator. Then, (A5'') is satisfied.\\[6pt]
The wellposedness of 2BSDEs will not alter under (A3'') and (A5''). First, the statement (\ref{er1}) remains true if we change a little of its expression:
$$
\mathbb{E}^\mathbb{P}_\tau\bigg[\int^1_\tau|\hat{a}^{1/2}_tZ_t|^2\bigg]
\leq\frac{1}{\gamma^2}e^{4\gamma||Y||_{\mathbb{D}^\infty_H}}(1+2\gamma (\overline{\alpha}+\overline{\beta}||Y||_{\mathbb{D}^\infty_H})),
$$
which yields that $Z$ is a $BMO(\mathcal{P}_H)$-martingale generator if $Y\in\mathbb{D}^\infty_H$. Lemma \ref{BMOP} and \ref{bmolemma2} ensure that the constants that we need for the proof of the representation theorem and the last step of the proof to the existence are uniform in $\mathbb{P}$. For the existence result, we have already explained that $V_t(\omega)$ in (\ref{v}) is well defined and all the properties still hold since (A3'') and (A5'') provide existence and uniqueness results as well as the estimates of solutions to quadratic BSDEs with the parameters $(\xi^{t, \omega}, \hat{F}^{t, \omega})$ under each $\mathbb{P}^t\in\mathcal{P}^t_H$.
\begin{remark}
If we assume in addition that $0\in\tilde{C}$, then $\overline{K}=\inf\{|r|: r\in a^{1/2}\tilde{C}\}=0$, so that the upper bound of $\hat{a}$ is not necessary. Both Theorem \ref{exptheo} and \ref{powtheo} can hold true under a larger class of probability measures $\hat{\mathcal{P}}_H$:
\begin{definition}\label{p4}
We denote by $\hat{\mathcal{P}}_H$ the collection which consists of all those $\mathbb{P}\in\overline{P}_S$ such that
$$\overline{a}^\mathbb{P}\leq \hat{a}_t\leq \underline{a}^\mathbb{P},\ {\rm tr}(\hat{a}_t^{-1})\leq\alpha_t,\ {\rm and}\ \hat{a}_t\in D_{F_t},\ \lambda\times\mathbb{P}-a.e.,$$
for some $\overline{a}^\mathbb{P}$, $\underline{a}^\mathbb{P}\in\mathbb{S}^{>0}_d$, a strictly positive $\alpha\in L^1([0,1])$ and each $(y, z)\in\mathbb{R}\times\mathbb{R}^d$. \\[6pt]Correspondingly, we denote by $\hat{\mathcal{P}}^t_H$ the collection of all those $\mathbb{P}^t\in\overline{P}^t_S$ such that
$$\overline{a}^{\mathbb{P}^t}\leq \hat{a}^t_s\leq \underline{a}^{\mathbb{P}^t},\ {\rm tr}((\hat{a}^t_s)^{-1})\leq\alpha_s\ {\rm and}\ \hat{a}^t_s\in D_{F_s},\ \lambda\times\mathbb{P}^t-a.e.,$$
for some $\overline{a}^{\mathbb{P}^t}$, $\underline{a}^{\mathbb{P}^t}\in\mathbb{S}^{>0}_d$, the same $\alpha$ as above and each $(y, z)\in\mathbb{R}\times\mathbb{R}^d$. 
\end{definition}
\noindent We can verify that this series $\{\mathcal{P}^t_H\}_{t\in [0, 1]}$ defined by Definition \ref{p4} is consistent and they ensure that (\ref{expgen}) and (\ref{powgen}) satisfy (A3'') and (A5''), respectively. \\[6pt]
However, if we consider the same problems under an even larger class of probability measures $\breve{\mathcal{P}}_H$:
\begin{definition}\label{p5}
We denote by $\breve{\mathcal{P}}_H$ the collection which consists of all $\mathbb{P}\in\overline{P}_S$ such that
$$\overline{a}^\mathbb{P}\leq \hat{a}_t\leq \underline{a}^\mathbb{P},\ \int^1_0 {\rm tr}(\hat{a}_t^{-1})dt\leq\overline{\alpha}\ {\rm and}\ \hat{a}_t\in D_{F_t},\ \lambda\times\mathbb{P}-a.e.,$$
for some $\overline{a}^\mathbb{P}$, $\underline{a}^\mathbb{P}\in\mathbb{S}^{>0}_d$, some strictly positive constant $\overline{\alpha}$ and each $(y, z)\in\mathbb{R}\times\mathbb{R}^d$.
\end{definition}
\noindent then the wellposedness of (\ref{exp2BSDE}) and (\ref{pow2BSDE}) will no longer hold true, since one is difficult to find a series of class $\{\breve{\mathcal{P}}^t_H\}_{t\in[0, 1]}$ consistent with $\breve{\mathcal{P}}_H$ defined by Definition \ref{p5}. In another word, once $\breve{\mathcal{P}}^t_H$ contains all the r.p.c.d. $\mathbb{P}^{t, \omega}$ of $\mathbb{P}\in\breve{\mathcal{P}}_H$, the second point in Definition \ref{consis} could not hold true. 
\end{remark}
\noindent\textbf{Acknowledgement} The author express special thanks to Prof. Hu, who provided both the initial inspiration for the work and useful suggestions.

\end{document}